\documentclass{m2an}

\usepackage{comment}
\usepackage{amsmath,tikz,latexsym}
\usepackage{amssymb}
\usepackage{graphicx}
\usepackage{color}
\usepackage{amsfonts}
\usepackage{amscd}
\usepackage{mathrsfs}
\usepackage{bm}
\usepackage{enumerate}
\usepackage{algorithmic, algorithm}

 \allowdisplaybreaks

\numberwithin{equation}{section}

    \newfont{\bg}{cmr10 scaled\magstep5}
    
    \newcommand{\bigzerou}{\smash{\lower1.7ex\hbox{\bg 0}}}

    \newcommand{\M}{ \mathcal{M}}

    \newcommand{\x}{{x}}
    \newcommand{\y}{{y}}
    \newcommand{\z}{{z}}
    \newcommand{\s}{\textbf{s}}
    \newcommand{\n}{\bm{n}}

     \newcommand{\vv}{\textbf{v}}

 \newcommand{\PP}{\textbf{P}}
\newcommand{\QQ}{\textbf{Q}}
 \newcommand{\LL}{\textbf{L}}

\begin{document}

\title{A Second-Order Nonlocal Approximation to Manifold Poisson Models with Neumann Boundary}
\thanks{This work was supported by NSFC grant 12071244, 11671005.}


\author{
Yajie Zhang}
\address{Department of Statistics and Mathematics, Zhongnan University of Economics and Law, Wuhan, China,
430000. Email: Z0005149@zuel.edu.cn}

\author{
Yanzun Meng}
\address{Department of Mathematical Sciences, Tsinghua University, Beijing, China,
100084. Email: yzmeng@mails.tsinghua.edu.cn}

\author{
Zuoqiang Shi}
\address{Corresponding Author, Department of Mathematical Sciences, Tsinghua University, Beijing, China,
100084. Email: zqshi@tsinghua.edu.cn}

\begin{abstract}
We propose a novel class of nonlocal models to approximate Poisson equations on manifolds with Neumann boundary conditions, where the manifolds are embedded in high-dimensional Euclidean spaces. Unlike existing approaches, our method optimizes truncation error by augmenting the model with a boundary-layer correction term involving second-order normal derivatives$-$expressed via the difference between interior and boundary Laplace-Beltrami operators$-$within the $2\delta$ interaction horizon. This work establishes three main results: (1) a geometrically consistent construction of the nonlocal model, (2) rigorous well-posedness analysis, and (3) proof of second-order convergence to the local Poisson model. Notably, our model achieves the optimal localization rate among known nonlocal approximations, including for manifolds in high-dimensional Euclidean spaces. Extensions of our model and several numerical tests are given after the main work.
 \end{abstract}


 \maketitle

\section{Introduction}

Partial differential equations on manifolds have garnered sustained research interest due to their broad applications across disciplines, including materials science \cite{CFP97} \cite{EE08}, fluid dynamics \cite{GT09} \cite{JL04}, biophysics \cite{BEM11}  \cite{ES10} \cite{NMWI11}, and more recently, machine learning \cite{belkin2003led} \cite{Coifman05geometricdiffusions}  \cite{LZ17}  \cite{MCL16}  \cite{reuter06dna} and image processing \cite{CLL15}    \cite{Gu04}  \cite{KLO17}  \cite{LWYGL14}   \cite{LDMM}   \cite{Peyre09}  \cite{Lui11}. The manifold Poisson equation has been particularly influential, serving both as a fundamental analytical tool for studying manifold geometry and as a practical model in these applications. While numerical analysis of Poisson models has seen significant advances, nonlocal approximations remain underexplored for manifolds compared to their Euclidean counterparts. Nonlocal models offer distinct advantages by eliminating explicit spatial derivatives, enabling innovative numerical schemes like the point integral method (PIM). This work addresses the critical need for high-accuracy nonlocal approximations of manifold Poisson equations to support emerging computational applications.

The development of nonlocal approximations for local models on Euclidean domains has matured considerably in recent years, with particular focus on
\begin{equation} \label{intro1}
\frac{1}{\delta^2} \int_{\Omega} (u_{\delta}(\x)-u_{\delta}(\y) ) R_\delta(\x,\y) d \y=f(\x), \qquad \x \in \Omega.
\end{equation}
Here $\Omega \subset \mathbb{R}^k$ is a bounded Euclid domain with smooth boundary,  $f \in H^2(\Omega)$, and $\delta$ represents the nonlocal interaction horizon. The nonlocal kernel is defined as $R_{\delta}(\x, \y) =C_{\delta} R \big( \frac{| \x-\y| ^2}{4 \delta^2} \big)$,
where $R \in C^2(\mathbb{R}^+) \cap L^1 [0, \infty)$ is a non-negative compactly support function, and  $C_{\delta}=\frac{1}{(4\pi\delta^2)^{k/2}}$ serves as the normalization constant. This formulation commonly appears in peridynamics literature \cite{Yunzhe4} \cite{Yunzhe8} \cite{Yunzhe13} \cite{Yunzhe27} \cite{Yunzhe31} \cite{Yunzhe32}.
 The truncation error between \eqref{intro1} and its local counterpart  $\Delta u=f$ exibits $\mathcal{O}(\delta^2)$ in interior regions but degrades to $\mathcal{O}(\delta^{-1})$ within the $2\delta$-boundary layer. Various modifications have been proposed to improve this boundary behavior: for Neumann conditions \cite{Yunzhe5} \cite{Yunzhe6} \cite{Yunzhe12} \cite{Yunzhe14}, other boundary types \cite{Yunzhe2} \cite{book-nonlocal} \cite{Du-SIAM} \cite{Yunzhe25}  \cite{ZD10},
 and alternative formulations using fractional kernels \cite{fractional1} or nonlocal gradients \cite{Stokes1}, though typically achieving only $\mathcal{O}(\delta)$ convergence.

Significant advances emerged with $\mathcal{O} (\delta^2)$ models for Neumann conditions in 1D \cite{Yunzhe} and 2D  \cite{Neumann_2nd_order}, later extended to Dirichlet problems via volumetric constraints \cite{Leehwi}. While \cite{YXB2} achieved second-order accuracy in higher dimensions, it requires first-order boundary data. Notably, \cite{YXB1} introduced polygonal-supported kernels that simplify numerical implementation.

However, nonlocal approximations on manifolds remain underdeveloped compared to Euclidean domains, limiting numerical options for manifold PDEs. Pioneering work in \cite{Base1} established a Neumann-type manifold model, followed by extensions to Dirichlet conditions \cite{Base2} \cite{WangTangJun}, interfaces \cite{Yjcms1}, and anisotropic cases \cite{Zqaniso}. The \cite{Base1} formulation employs
\begin{equation}  \label{base}
-\int_{\mathcal{M}} \Delta_{\mathcal{M}} u(\y) \bar{R}_{\delta}(\x,\y) d \mu_{\y} \approx \int_\M \frac{1}{\delta^2} R_{\delta}(\x,\y) (u(\x)-u(\y)) d \mu_{\y} -2 \int_{\partial \mathcal{M}} \bar{R}_{\delta}(\x,\y) \frac{\partial u}{\partial \n} (\y) d \tau_\y,
\end{equation}
Here $\M$ is a compact smooth $m$ dimensional manifold embedded in $\mathbb{R}^d$, with boundary $\partial \M$ a smooth $(m-1)$ dimensional submanifold. 
$\Delta_{\M}$ is the Laplace-Beltrami operator on $\M$, with $\nabla_{\M}$ the surface gradient(see \cite{Base1} $pp. 2-3$ for definitions). $u \in H^4(\M)$,
$\frac{\partial u}{\partial \n} =\nabla_{\M} u \cdot \n$ is the normal derivative, and $d \mu_\y$ and $d \tau_\y$ are the volume forms of $\M$ and $\partial \M$, respectively.
The kernel function $R_{\delta}$ are defined as $ R_{\delta}(\x, \y) =C_{\delta} R \big( \frac{| \x-\y| ^2}{4 \delta^2} \big), \mbox{with } C_{\delta}=(4\pi\delta^2)^{-m/2},
\mbox{ and } R \in C^2(\mathbb{R}^+) \cap L^1 [0, \infty)$ is a non-negative, compactly supported function.
 $ \bar{R}_{\delta}(\x, \y) =C_{\delta} \bar{R} \big( \frac{| \x-\y| ^2}{4 \delta^2} \big)$,  $\bar{R}(r)=\int_r^{+\infty} R(s) ds$. 
A key approximation $2\delta^2 \nabla^{\x}_{\M} \bar{R}_{\delta}(\x,\y) \approx (\y-\x) R_{\delta}(\x,\y)$ holds for $\x,\y \in \M$.

Theorem 4 of \cite{Base1} provides a preliminary error analysis of \eqref{base}. The Neumann boundary model is derived by omitting the last term in \eqref{base} and balancing the resulting equation, achieving $\mathcal{O}(\delta)$ convergence to the local counterpart. A similar rate was obtained for Dirichlet boundaries in \cite{LSS}, where the boundary condition was approximated via a Robin-type condition $\frac{\partial u}{\partial \n} \approx \lambda u$.

To achieve $\mathcal{O} (\delta^2)$ convergence$-$critical for high-accuracy applications$-$we proposed in  \cite{Base3} \cite{Base4} a modified model incorporating the second-order normal derivative of $u$ along 
$\partial \M$. This term corrects the dominant boundary-layer error in \eqref{base}, enabling second-order accuracy:
 \begin{equation}  \label{b0016}
 \begin{split}
 \!\!\!\!\!\!\!\!\! -\int_{\M} \Delta_{\M} u(\y) \bar{R}_{\delta}(\x,\y) d \mu_\y = \ & \frac{1}{\delta^2} \int_\M (u(\x)-u(\y)) R_{\delta}(\x,\y)  d \mu_\y -2 \int_{\partial \M}  \frac{\partial u}{\partial \n} (\y) \bar{R}_{\delta}(\x,\y) d \tau_\y \\
 & - \int_{\partial \M} (\x-\y) \cdot  \n(\y)    \nabla_{\M}^2 u (\n,\n) (\y)     \bar{R}_{\delta} (\x, \y) d \tau_\y
 -r_{in}(\x)  , \ 
 \x \in \M,
 \end{split}
\end{equation}
the term $r_{in}$, which is $\mathcal{O}(\delta^2)$, was fully characterized in Theorem 3.1 of \cite{Base4}. 
Building on \eqref{b0016}, we developed in \cite{Base3} a nonlocal Poisson model with $\mathcal{O}(\delta^2)$ convergence for Dirichlet boundary conditions. In this case, the second-order normal derivative $ \nabla_{\M}^2 u (\n,\n)$ can be expressed in terms of $\frac{\partial u}{\partial \n}$, which itself admits a boundary-layer convolution approximation.

The absence of explicit second-order normal derivative information for Neumann boundary conditions presents a fundamental obstacle in constructing nonlocal models from \eqref{b0016}. While \cite{Neumann_2nd_order} resolved this for 2D Euclidean domains by expressing the derivative through level-curve integrals of $u-$yielding an $\mathcal{O}(\delta^2)$-accurate well-posed model$-$their technique fails to generalize to higher-dimensional manifolds.
 Our work addresses the fundamental challenge of handling second-order normal derivatives in Neumann boundary problems through two key innovations: first, by representing the derivative via the difference between interior and boundary Laplace-Beltrami operators; second, by approximating the boundary Laplacian through a carefully constructed augmented function that preserves model coercivity. The homogeneous Neumann condition (vanishing first-order derivative) allows us to rigorously reformulate \eqref{b0016}, yielding a new nonlocal formulation with optimal convergence properties.

 This work makes three primary contributions: the construction of a new nonlocal model, proof of its well-posedness, and demonstration of its second-order convergence rate to the local counterpart. Notably, this represents the first nonlocal Poisson model achieving $\mathcal{O}(\delta^2)$  convergence for Neumann boundary conditions in dimensions $d \geq 3$, offering significant improvements in numerical efficiency. The model's compatibility with meshless schemes like the Point Integral Method (PIM) \cite{LSS} provides distinct advantages for manifold problems, particularly where traditional finite-element methods (FEM) become impractical due to high-dimensional mesh generation complexities. Beyond immediate applications, this work establishes a foundational approach for handling nonlocal problems on manifolds.
 
The paper is organized as follows: Section 2 develops the formulation of our nonlocal Poisson model. Section 3 presents the well-posedness analysis and Sobolev regularity improvements. The $\mathcal{O}(\delta^2)$ convergence rate to the local counterpart is established in Section 4. Section 5 discusses model generalizations, while Section 6 presents numerical examples to validate such rate of convergence. We conclude with a summary and future research directions in Section 7. For readability, several technical proofs appear in the appendix.
   
\section{Construction of Nonlocal Model}
This work develops a second-order nonlocal approximation for the following Neumann problem defined on a smooth manifold:
\begin{equation}  \label{bg01}
\begin{cases}
-\Delta_{\M} u(\x) =f(\x) & \x \in \M, \\
\frac{\partial u}{\partial \n}(\x)=0  & \x \in \partial \M, \\
\int_{\M} u(\x) d \mu_{\x}=0.
\end{cases}
\end{equation}
where $\M$ is an open, bounded, smooth $m$-dimensional manifold embedded in $\mathbb{R}^d$. Its boundary, $\partial \M$, is a smooth $(m-1)$-dimensional submanifold. Here, $\Delta_{\M}$ denotes the Laplace-Beltrami operator on $\M$, and the source term $f$ belongs to $H^3(\M)$ and satisfies $\int_{\M} f (\x) d \mu_{\x} =0$. Classical elliptic regularity theory guarantees that \eqref{bg01} admits a unique solution $u \in H^5(\M)$.

To construct the nonlocal model, we first define a kernel function of the form:
$$R_{\delta}(\x, \y) =C_{\delta} R \big( \frac{| \x-\y| ^2}{4 \delta^2} \big).$$
where $\delta>0$ is the interaction horizon and $C_{\delta}=(4\pi\delta^2)^{-m/2}$ is a normalization constant.
To facilitate an efficient error analysis, we impose the following conditions on the radial profile $R(r)$:
\begin{enumerate}
\item Smoothness: $R \in W^{2, +\infty}[0, +\infty)$, i.e., for any $r \geq 0$ we have $\big| D^2 R(r) \big| \leq C$;
\item Nonnegativity: $R(r) \geq 0$ for any $r \geq 0$;
\item Compact support: $R(r)=0$ for any $r>1$;
\item Nondegenearcy: $\exists \ \delta_0 >0$ so that $R(r) \geq \delta_0 >0$ for $0 \leq r \leq 1/2$.
\end{enumerate}
These assumptions are general and are satisfied by most standard kernels in the literature. The smoothness condition prevents singular behavior near the origin, while the non-degeneracy condition ensures the kernel retains a substantial support region. Although relaxing the compact support requirement leads to more challenging theoretical problems, we retain it here to streamline the error estimates.

A canonical example of a kernel satisfying these requirements is the cosine kernel:
\begin{equation}
R(r)=
\begin{cases}
\frac{1}{2}(1+ \cos \pi r ), & 0 \leq r \leq 1, \\
0, & r>1.
\end{cases}
\end{equation}
Furthermore, we define three auxiliary functions generated from $R$ via successive integration:
\begin{equation} \label{kernelori}
{R}(r)=\int_r^{+\infty} \underline{R}(s) ds, \quad \bar{R}(r)=\int_r^{+\infty} R(s) ds, \quad \overset{=}{R}(r)=\int_r^{+\infty} \bar{R}(s) ds.
\end{equation}
The corresponding nonlocal kernels are then given by:
\begin{equation} \label{kernel}
\underline{R}_{\delta}(\x, \y) =C_{\delta} \underline{R} \big( \frac{| \x-\y| ^2}{4 \delta^2} \big), \quad 
\bar{R}_{\delta}(\x, \y) =C_{\delta} \bar{R} \big( \frac{| \x-\y| ^2}{4 \delta^2} \big) , \quad  \overset{=}{R}_{\delta}(\x, \y) =C_{\delta} \overset{=}{R} \big( \frac{| \x-\y| ^2}{4 \delta^2} \big).
\end{equation}
A direct verification confirms that these derived kernels also satisfy the original assumptions, with only minor modifications required for the smoothness and non-degeneracy conditions.

We now derive a nonlocal approximation for the system \eqref{bg01}. Let $ u \in H^5(\M) $ be its unique solution. Incorporating the equality $ f = -\Delta_{\M} u$ and the boundary condition $ \frac{\partial u}{\partial \n} = 0 $ into \eqref{b0016}, we obtain:
 \begin{equation}  \label{b0017}
 \begin{split}
  \int_{\M} f(\y) \bar{R}_{\delta}(\x,\y) d \mu_\y  = \frac{1}{\delta^2} \int_\M (u(\x)-u(\y)) R_{\delta}(\x,\y)  d \mu_\y 
  -  \int_{\partial \M} (\x-\y) \cdot  \n(\y)   \nabla_{\M}^2 u (\n,\n) (\y)     \bar{R}_{\delta} (\x, \y) d \tau_\y -r_1(\x) ,  \ \x \in \M;
 \end{split}
\end{equation}
where we have relabeled the approximation error as $r_1$ (replacing $r_{in}$ from \eqref{b0016}) for notational consistency.

Next, we handle the second-order term. The following identity from differential geometry holds:
\begin{equation} \label{mm10}
\nabla_{\M}^2 u (\n,\n) (\y)   =\Delta_\M u(\y) - \Delta_{\partial \M} u(\y)
\end{equation}
under the constraint $\frac{\partial u}{\partial \n} = 0$. The proof of \eqref{mm10} is given in Appendix \ref{appen0}.
Substituting $\Delta_\M=-f$ into \eqref{mm10} and applying the result to \eqref{b0017} yields
 \begin{equation}  \label{b0019}
 \begin{split}
& \frac{1}{\delta^2} \int_\M (u(\x)-u(\y))  R_{\delta}(\x,\y)  d \mu_\y  +   \int_{\partial \M} \Delta_{ \partial \M} u(\y)   (\x-\y) \cdot  \n(\y)  \bar{R}_{\delta} (\x, \y)      d \tau_\y  \\
& =  \int_{\M} f(\y) \bar{R}_{\delta}(\x,\y) d \mu_\y - \int_{\partial \M} (\x-\y) \cdot  \n(\y) f(\y)     \bar{R}_{\delta} (\x, \y) d \tau_\y+r_1(\x), \ \x \in \M.
 \end{split}
\end{equation}
For brevity, define
\begin{equation} \label{fdelta}
{f}_{\delta}(\x)= \int_{\M} f(\y) \bar{R}_{\delta}(\x,\y) d \mu_\y - \int_{\partial \M} (\x-\y) \cdot  \n(\y) f(\y)     \bar{R}_{\delta} (\x, \y) d \tau_\y, \quad \x \in \M;
\end{equation}
\begin{equation}
\zeta'_{\delta} (\x,\y)=-(\x-\y) \cdot  \n(\y)     \bar{R}_{\delta} (\x, \y),   \quad \x \in \M, \ \y \in \partial \M.
\end{equation}
Neglecting the truncation error $r_1$ and replacing $u$ with $u_{\delta}$ in \eqref{b0019}, we obtain the following preliminary model:
  \begin{equation}  \label{preli}
 \begin{split}
 \frac{1}{\delta^2} \int_\M (u_{\delta}(\x)-u_{\delta}(\y))  R_{\delta}(\x,\y)  d \mu_\y  -   \int_{\partial \M} \Delta_{ \partial \M} u_{\delta}(\y)   \zeta'_{\delta}(\x,\y)   d \tau_\y   =  {f}_{\delta}(\x), \quad \x \in \M.
 \end{split}
\end{equation}
This model couples the boundary Laplacian $\Delta_{ \partial \M} u_{\delta}(\y)$ to the geometric factor 
$\zeta'_{\delta}(\x,\y)$.
However, in the associated energy, this term becomes non-sign-definite after integration by parts, leading to potential instabilities in numerical schemes and violating the physical principle of energy dissipation for diffusion processes. Specifically, the bilinear form from this term fails to satisfy non-negativity, rendering the energy functional potentially unbounded below.
To ensure energy non-negativity while preserving geometric boundary interactions, we modify the preliminary model 
 \eqref{preli} in four steps:
\begin{enumerate}
\item  \textbf{Boundary Laplacian.}  Introduce an auxiliary boundary potential $\hat{u}'_{\delta}(\y)$, defined as a kernel-smoothed projection of the bulk field $u_{\delta}$ onto the boundary geometry:
\begin{equation} \label{preli5}
\hat{u}'_{\delta}(\y)=\frac{\int_{\M} u_{\delta}(\x) \zeta'_{\delta}(\x,\y) d \mu_{\x}  } {\int_{\M}  \zeta'_{\delta}(\x,\y) d \mu_{\x}   }
\end{equation}
This variable $\hat{u}'_{\delta}$ captures the local imbalance of $u_{\delta}$ relative to the boundary orientation. 
Replace the original boundary term  $\Delta_{ \partial \M} u_{\delta}(\y)$ in \eqref{preli} with
$ \Delta_{ \partial \M} \hat{u}'_{\delta}(\y)$. The boundary term's contribution to the total energy becomes:
\begin{equation} \label{preli2}
-\int_{\M} u_{\delta} (\x)   \int_{\partial \M} \Delta_{ \partial \M} \hat{u}'_{\delta}(\y)   \zeta'_{\delta}(\x,\y)   d \tau_\y  d \mu_{\x} =- \int_{\partial \M} \Delta_{ \partial \M} \hat{u}'_{\delta}(\y) \ \hat{u}'_{\delta}(\y)  \ \bigg( \int_{\M}   \zeta'_{\delta}(\x,\y) d \mu_{\x}   \bigg) \ d \tau_{\y}
\end{equation}
\item \textbf{Geometric Factor. }
The energy functional \eqref{preli2} remains problematic due to the presence of $\int_{\M}   \zeta'_{\delta}(\x,\y) d \mu_{\x}$. To address this, we modify the geometric factor $\zeta'_{\delta}$ in \eqref{preli}. Note that
\begin{equation} 
\int_{\M} \zeta'_{\delta}(\x,\y) \, d\mu_{\x} \approx 2\delta^2 \int_{\M} \nabla^{\x}_{\M} \overset{=}{R}_{\delta}(\x,\y) \cdot \n(\y) \, d\mu_{\x} \approx 2\delta^2 \int_{\partial \M} \overset{=}{R}_{\delta}(\x,\y) \, d\tau_{\x} \approx \delta C_R,
\end{equation}
where $ C_R = \pi^{-m/2} \int_{\mathbb{R}^{m-1}} \overset{=}{R}(|\x|^2) \, d\x $.
Replace $\zeta'_{\delta}(\x,\y)$ in \eqref{preli} with
\begin{equation} \label{preli4}
\zeta'_{\delta} (\x,\y) \frac{ \delta C_R}{ \int_{\M} \zeta'_{\delta}(\x,\y) d \mu_{\x} } .
\end{equation}
The modified boundary term's energy contribution becomes:
\begin{equation} \label{preli3}
-\int_{\M} u_{\delta} (\x)   \int_{\partial \M} \Delta_{ \partial \M} \hat{u}'_{\delta}(\y) \frac{ \delta C_R \   \zeta'_{\delta}(\x,\y)  }{ \int_{\M} \zeta'_{\delta}(\z,\y) d\mu_{\z}}  d \tau_\y  d \mu_{\x}= \delta C_R \left \| \nabla_{\partial \M} \hat{u}'_{\delta} \right \|^2_{L^2(\partial \M)}.
\end{equation}
which is non-negative.
\item \textbf{Simplification.}
To simplify the model, replace $\zeta'_{\delta}$ with
\begin{equation} \label{zeta1}
\begin{split}
\zeta_{\delta}(\x,\y)=2\delta^2 \n(\y) \cdot \nabla^{\x}_{\M} \overset{=}{R}_{\delta} (\x, \y)   .
\end{split} 
\end{equation}
Note that $\zeta_{\delta} =\zeta'_{\delta}$ when $\M$ is a Euclid domain. This choice simplifies integration over $\M$. Define the updated auxiliary variable:
\begin{equation} \label{preli6}
\hat{u}_{\delta}(\y)=\frac{\int_{\M} u_{\delta}(\x) \zeta_{\delta}(\x,\y) d \mu_{\x}  } {\int_{\M}  \zeta_{\delta}(\x,\y) d \mu_{\x}   }.
\end{equation}
The final boundary integral term replacing the original in \eqref{preli} is
\begin{equation}  \label{preli7}
- \int_{\partial \M} \Delta_{ \partial \M} \hat{u}_{\delta}(\y)   \zeta_{\delta}(\x,\y) \frac{ \delta C_R}{ \int_{\M} \zeta_{\delta}(\z,\y) d\mu_{\z}}  d \tau_\y.
 \end{equation}
 
 \item \textbf{Average-Zero Constraint.}
 Using integration by parts, both  \eqref{preli7} and the nonlocal diffusion term in \eqref{preli} are average-zero over $\M$. To ensure solution existence, subtract the constant:
 \begin{equation} \label{preli8}
 \tilde{f}_{\delta}=\frac{1}{\mbox{Vol}(\M)} \int_{\M} f_{\delta}(\x) d\mu_{\x}
 \end{equation}
from the right-hand side of \eqref{preli} to enforce average-zero conditions on both sides.
 
\end{enumerate}

In conclusion, we propose the following nonlocal model:
 \begin{equation}  \label{c01}
 \begin{cases}
\displaystyle  \frac{1}{\delta^2}  \int_\M (u_{\delta}(\x)-u_{\delta}(\y)) R_{\delta}(\x,\y)  d \mu_\y  -  \delta C_R \int_{\partial \M} \Delta_{ \partial \M} \hat{u}_{\delta}(\y)   \zeta_{\delta}(\x,\y)  \frac{1}{\hat{\omega}_{\delta}(\y)} d \tau_\y 
= f_{\delta}(\x)-\tilde{f}_{\delta}  \qquad & \x \in \M; \\
 \displaystyle  \hat{u}_{\delta}(\y)- \frac{1}{\hat{\omega}_{\delta}(\y)} \int_{\M} u_{\delta}(\x) \zeta_{\delta}(\x,\y) d \mu_{\x}=0 , \qquad & \y \in \partial \M; \\
\displaystyle \int_{\M} u_{\delta}(\x) d \mu_{\x}=0,
 \end{cases}
\end{equation}
where $\zeta_{\delta}, \tilde{f}_{\delta}$ are defined in \eqref{zeta1} \eqref{preli8} respectively, and:
\begin{equation} \label{defome1}
 \hat{\omega}_{\delta}(\y)=\int_{\M}    \zeta_{\delta}(\x,\y)    d \mu_\x , \ \forall \ \y \in \partial \M; \qquad 
C_R=  \pi^{-\frac{m}{2}} \int_{\mathbb{R}^{m-1}} \overset{=}{R}(|\x|^2) d\x.
\end{equation}

Now we are ready to state our main theorem.
\begin{thrm} \label{theorem1}
Let $\mathcal{M}$ be an $m$-dimensional $C^4$-smooth manifold embedded in $\mathbb{R}^n$ with $C^4$-smooth $(m-1)$-dimensional boundary $\partial \mathcal{M}$. Then:
\begin{enumerate}
\item (Well-Posedness) For any $f \in H^1(\mathcal{M})$, there exists a unique solution $u_{\delta} \in H^1(\mathcal{M})$ to the nonlocal model \eqref{c01} satisfying:
   \begin{equation}
  \left \| u_{\delta} \right \|_{H^1(\mathcal{M})} \leq C \left \| f \right \|_{H^1(\mathcal{M})},
   \end{equation}
   where $C > 0$ depends only on the geometry of $\mathcal{M}$ and $\partial \mathcal{M}$.

\item (Second-Order Convergence) If $f \in H^3(\mathcal{M})$ and $u$ solves the local problem \eqref{bg01}, then:
    \begin{equation}
   \left \| u - u_{\delta} \right \|_{H^1(\mathcal{M})} \leq C \delta^2 \left \| f \right \|_{H^3(\mathcal{M})},
      \end{equation}
   with $C > 0$ again depending only on geometric properties.
\end{enumerate}
\end{thrm}

\section{Well-Posedness} \label{psd}
In this section, we prove part (i) of Theorem 2.1. We begin by defining the necessary function spaces and bilinear forms.

Let $L_0^2(\M)$ denote the Hilbert space of $L^2$ functions on the manifold $\M$ with zero mean. We equip it with the standard $L^2$ inner product:
\begin{equation}
( p_{\delta}, q_{\delta})= \int_{\M} p_{\delta}(\x) q_{\delta}(\x) d\mu_{\x}.
\end{equation}
The core of our analysis is the following bilinear form on $L_0^2(\M)$:
\begin{equation} \label{bilinear2}
 \begin{split}
 B_{\delta}[p_{\delta}, q_{\delta}]= 
\frac{1}{\delta^2} \int_\M q_{\delta}(\x) \int_\M  (p_{\delta}(\x)-p_{\delta}(\y))  R_{\delta}(\x,\y)  d \mu_\y d \mu_\x  
 - \delta C_R \int_\M q_{\delta}(\x)   \int_{\partial \M} \Delta_{ \partial \M} \hat{p}_{\delta}(\y)   \zeta_{\delta}(\x,\y)  \frac{1}{\hat{\omega}_{\delta}(\y)} d \tau_\y d \mu_\x,
 \end{split}
\end{equation}
where
 \begin{equation} 
 \hat{p}_{\delta} (\y) =\frac{ 1 }{ \hat{\omega}_{\delta}(\y) } \int_{\M} p_{\delta}(\z)  \zeta_{\delta}(\z,\y) d \mu_{\z}.
\end{equation}
The well-posedness of the associated variational problem is established by verifying the conditions of the Lax-Milgram theorem, which we now do in the following lemma.

\begin{lmm} \label{Lax}
The bilinear form $B_{\delta}[\cdot, \cdot]$ is bounded and coercive on $L^2_0(\M)$. Specifically, there exist constants $C_{\delta}>0$ and $C>0$, where $C_{\delta}$ depends on $\delta$ and the manifold $\M$, and $C$ depend only on the geometry of $\M$, such that for all $p_{\delta}, q_{\delta} \in L_0^2(\M) $:
\begin{enumerate} 
\item Boundedness:
\begin{equation} \label{lax2}
|B_{\delta}[p_{\delta}, q_{\delta}]| \leq C_{\delta} \left \| p_{\delta} \right \|_{L^2(\M)} \left \| q_{\delta} \right \|_{L^2(\M)}, \qquad \forall \ p_{\delta}, q_{\delta} \in L_0^2(\M);
\end{equation}
\item Coercivity:
\begin{equation} \label{lax3}
B_{\delta}[p_{\delta}, p_{\delta}]  \geq C \left \| p_{\delta} \right \|^2_{L^2(\M)}, \qquad \forall \ p_{\delta} \in L_0^2(\M);
\end{equation}
\end{enumerate}
\end{lmm}

\begin{proof}[Proof of Lemma \ref{Lax}]
The definition of $B_{\delta}$ presented in \eqref{bilinear2} yields the following estimate:
\begin{equation} \label{bilinear4}
\begin{split}
 B_{\delta}[p_{\delta}, q_{\delta}] \leq & C_{\delta} ( \int_\M \int_\M |q_{\delta}(\x)|  (| p_{\delta}(\x)|+| p_{\delta}(\y)|) d \mu_\y d \mu_\x 
 +  C_{\delta}  \int_\M |q_{\delta}(\x)| \int_{\partial \M} |\hat{p}_{\delta} (\y)  \Delta_{ \partial \M}  (\zeta_{\delta}(\x,\y)  \frac{1}{\hat{\omega}_{\delta}(\y)}) | d \tau_{\y} d\mu_{\x}  \\
\leq  & C_{\delta} (\left \| p_{\delta} \right \|_{L^2(\M)} \left \| q_{\delta} \right \|_{L^2(\M)}
 + \left \| p_{\delta} \right \|_{L^1(\M)} \left \| q_{\delta} \right \|_{L^1(\M)}  
 + \left \| \hat{p}_{\delta} \right \|_{L^1(\M)} \left \| q_{\delta} \right \|_{L^1(\M)} )
 \leq C_{\delta}  \left \| p_{\delta} \right \|_{L^2(\M)} \left \| q_{\delta} \right \|_{L^2(\M)}.
\end{split}
\end{equation}
Now our only goal is to prove \eqref{lax3}. Setting $q_{\delta}=p_{\delta}$, we simplify \eqref{bilinear2} into
\begin{equation} \label{bilinear5}
\begin{split}
B_{\delta}[p_{\delta}, p_{\delta}] = &
\frac{1}{\delta^2} \int_\M \int_\M p_{\delta}(\x) (p_{\delta}(\x)-p_{\delta}(\y))  R_{\delta}(\x,\y)  d \mu_\y d \mu_\x  
 - \delta C_R \int_{\partial \M} \Delta_{ \partial \M} \hat{p}_{\delta}(\y)    \int_{ \M}  p_{\delta}(\x)  \zeta_{\delta}(\x,\y)  \frac{1}{\hat{\omega}_{\delta}(\y)} d \mu_\x d \tau_\y \\
 = & \frac{1}{2\delta^2} \int_\M \int_\M (p_{\delta}(\x) -p_{\delta}(\y) )^2  R_{\delta}(\x,\y)  d \mu_\y d \mu_\x  
 + \delta C_R \int_{\partial \M}   \nabla_{ \partial \M} \hat{p}_{\delta}(\y)  \cdot \nabla_{ \partial \M} \hat{p}_{\delta}(\y) d \tau_{\y}.
\end{split}
\end{equation}
Both terms are manifestly non-negative.

To prove coercivity, we define a weighted-average function:
\begin{equation} \label{tildep}
\tilde{p}_{\delta} (\x) = \frac{1}{\omega_1(\x)} \int_{\M} p_{\delta} (\y) R_{\delta}(\x,\y)  d \mu_\y,  \qquad
\mbox{where} \ \omega_1(\x)= \int_{\M} R_{\delta}(\x,\y)  d \mu_\y, 
\qquad \forall \ \x \in \M.
 \end{equation}
The Poincar\'{e} inequality on $\M$ implies:
\begin{equation} \label{bilinear7}
\begin{split}
 \left \| p_{\delta} \right \|^2_{L^2(\M)}
 \leq C( \left \| p_{\delta}-\tilde{p}_{\delta} \right \|^2_{L^2(\M)} + \left \| \tilde{p}_{\delta} \right \|^2_{L^2(\M)})
 \leq C( \left \| p_{\delta}-\tilde{p}_{\delta} \right \|^2_{L^2(\M)} +  ( \int_{\M} \tilde{p}_{\delta}(\x)  d \mu_{\x})^2+ \left \| \nabla_{\M} \tilde{p}_{\delta} \right \|^2_{L^2(\M)}  ).
 \end{split}
 \end{equation}
 We control the right-hand side term-by-term. A standard nonlocal-to-local estimate (cf. \cite{Yjcms1}, Eq. (5.17)) yields
\begin{equation} \label{bilinear6}
 \left \| \nabla_{\M} \tilde{p}_{\delta} \right \|_{L^2(\M)}^2 \leq  \frac{C}{2\delta^2} \int_\M \int_\M (p_{\delta}(\x) -p_{\delta}(\y) )^2  R_{\delta}(\x,\y)  d \mu_\y d \mu_\x \leq B_{\delta}[p_{\delta}, p_{\delta}];
\end{equation}
in addition, noticing that $\int_{\M} p_{\delta}(\x) d \mu_{\x}=0$, we have
\begin{equation} \label{bilinear8}
\begin{split}
 (\int_{\M} \tilde{p}_{\delta}(\x) & d \mu_{\x} )^2= (\int_{\M} \tilde{p}_{\delta}(\x) -p_{\delta}(\x) d \mu_{\x} )^2
=  \big( \int_{\M}   \int_{\M} \frac{1}{\omega_1(\x)} (p_{\delta} (\y)-p_{\delta}(\x)) R_{\delta}(\x,\y)  d \mu_\y d \mu_{\x}  \big)^2 \\
& \leq C \int_{\M}  \int_{\M} (p_{\delta} (\y)  -p_{\delta}(\x))^2 R_{\delta}(\x,\y)  d \mu_\y d{\mu}_{\x} \int_{\M} \int_{\M}   \frac{1}{\omega^2_1(\x)}  R_{\delta}(\x,\y)  d \mu_\y d \mu_{\x}  
\leq C \delta^2 B_{\delta}[p_{\delta}, p_{\delta}];
\end{split}
\end{equation}
and similarly
\begin{equation} \label{bilinear9}
\begin{split}
  \left \| {p}_{\delta}-  \tilde{p}_{\delta} \right \|^2_{L^2(\M)} &
 =   \int_{\M} \big(   \int_{\M} \frac{1}{\omega_1(\x)} (p_{\delta} (\y)-p_{\delta}(\x)) {R}_{\delta}(\x,\y)  d \mu_\y \big)^2 d \mu_{\x}   \\
& \leq C \int_{\M}  (\int_{\M} (p_{\delta} (\y)-p_{\delta}(\x))^2 {R}_{\delta}(\x,\y)  d \mu_\y ) ( \int_{\M}    \frac{1}{\omega^2_1(\x)}  {R}_{\delta}(\x,\y)  d \mu_\y )d \mu_{\x}    \leq C \delta^2 B_{\delta}[p_{\delta}, p_{\delta}], 
\end{split}
\end{equation}
we thereby establish \eqref{lax3} by combining \eqref{bilinear7}$-$\eqref{bilinear9}.
\end{proof}

With Lemma \ref{Lax} established, the bilinear form $B_{\delta}[ \cdot, \cdot]$ satisfies the hypotheses of the Lax-Milgram theorem.

Furthermore, the definition of $f_{\delta}$ in \eqref{fdelta} yields the estimate
\begin{equation} \label{higher15}
\begin{split}
& \left \| f_{\delta} \right \|^2_{L^2(\M)} 
=    \int_{\M} ( \int_{\M} f(\y) \bar{R}_{\delta}(\x,\y) d \mu_\y - \int_{\partial \M} (\x-\y) \cdot  \n(\y)   f(\y)      \bar{R}_{\delta} (\x, \y) d \tau_\y )^2 d \mu_{\x}  \\
\leq & C   \int_{ \M}  (\int_{\M} f^2(\y) \bar{R}_{\delta}(\x,\y) d \mu_\y) (\int_{\M} \bar{R}_{\delta}(\x,\y) d \mu_\y) d \mu_{\x} + 
C \int_{ \M}  (\int_{\partial \M} \delta^2 f^2(\y) \bar{R}_{\delta}(\x,\y) d \tau_\y) (\int_{\partial \M} \bar{R}_{\delta}(\x,\y) d \tau_\y) d \mu_{\x}  \\
\leq & C   ( \int_{\M}  \int_{\M} f^2(\y) \bar{R}_{\delta}(\x,\y) d \mu_\y d \mu_{\x}  + 
 \delta \int_{ \M}  \int_{\partial \M} f^2(\y) \bar{R}_{\delta}(\x,\y) d \tau_\y d \mu_{\x}  \big) 
\leq  C (\left \| f \right \|^2_{L^2(\M)} + \delta \left \| f \right \|^2_{L^2(\partial \M)} )
\leq C  \left \| f \right \|^2_{H^1(\M)}.
\end{split}
\end{equation}
Consequently, its spatial average, $\tilde{f}_{\delta}$ (defined in \eqref{preli8}), is also bounded by $C  \left \| f \right \|_{H^1(\M)}$ as well. According to Lax-Milgram theorem and Lemma \ref{Lax}, the problem of finding $u_{\delta} \in L_0^2(\M)$ such that
\begin{equation} \label{bilinear1}
B_{\delta}[u_{\delta}, q_{\delta}]=(f_{\delta}-\tilde{f}_{\delta}, q_\delta) \qquad \forall \ q_{\delta} \in L_0^2(\M)
\end{equation}
admits a unique solution. This solution satisfies the bound
\begin{equation} \label{LaxMilgramfinal}
\left \| u_{\delta} \right \|_{L^2(\M)} \leq C \left \| f_{\delta}-\tilde{f}_{\delta} \right \|_{L^2(\M)}  \leq C \left \| f \right \|_{H^1(\M)}.
\end{equation}
We now demonstrate the equivalence between the weak formulation \eqref{bilinear1} and the strong form \eqref{c01}.
\begin{enumerate} 
\item \textbf{Strong implies Weak:} Suppose $u_{\delta}$ is a solution to \eqref{c01}. Taking the inner product of the first equation in \eqref{c01} with an arbitrary test function $q_{\delta} \in L_0^2(\M)$ and integrating over $\M$ immediately yields the weak form \eqref{bilinear1}.
\item \textbf{Weak implies Strong:} 
Now, suppose $u_{\delta} \in L_0^2(\M)$ satisfies the weak formulation \eqref{bilinear1} for all $q_{\delta} \in L_0^2(\M)$. A key observation is that $B_{\delta}[u_{\delta},1]=(f_{\delta}-\tilde{f}_{\delta},1)=0$ by construction.
This allows us to extend the validity of \eqref{bilinear1} from $q_{\delta} \in L_0^2(\M)$ to the entire space $L^2(\M)$. By the fundamental lemma of the calculus of variations, this implies that the first equation in \eqref{c01} holds almost everywhere in $\M$. The second and third equations in \eqref{c01} are enforced by the definition of the function space $L_0^2(\M)$ and the boundary condition built into the formulation.
 \end{enumerate}

 To complete the proof of Theorem \ref{theorem1}(i), it remains to establish the higher regularity $u_{\delta} \in H^1(\mathcal{M})$. This will follow from the subsequent lemma.
\begin{lmm} \label{higher1}
Let $r_{\delta} \in H^1(\M)$, and suppose $p_{\delta} \in L^2(\M)$ satisfies the bilinear form
\begin{equation} \label{higher2}
B_{\delta}[p_{\delta}, q_{\delta}]=(r_{\delta}, q_\delta) \qquad \forall \ q_{\delta} \in L^2(\M),
\end{equation}
where 
\begin{equation} \label{higher3}
\!\!\!\!\!\!\!\!\!\!\!\! (r_{\delta}, q_\delta) \leq C \left \| r_{1\delta} \right \|_{H^k(\M)} (\left \| q_{\delta} \right \|_{H^1(\M)} + \left \| \bar{q}_{\delta} \right \|_{H^1(\M)}  + \left \| \overset{=}{q}_{\delta} \right \|_{H^1(\M)} ),
\ \forall \ q_{\delta} \in H^1(\M) 
\end{equation}
for some $k \in \mathbb{N}^+$ and $r_{1\delta} \in H^m(\M)$. Here 
\begin{equation} \label{omega2}
\bar{q}_{\delta} (\x) =\frac{ 1 }{ \omega_2(\x) } \int_{\M} q_{\delta} (\y) \bar{R}_{\delta}(\x,\y) d \mu_{\x}, 
\qquad \omega_2(\x)=  \int_{\M} \bar{R}_{\delta}(\x,\y) d \mu_{\x},
\end{equation}
\begin{equation} \label{omega3}
\overset{=}{q}_{\delta} (\x) =\frac{ 1 }{ \omega_3(\x) } \int_{\M} q_{\delta} (\y) \overset{=}{R}_{\delta}(\x,\y) d \mu_{\x}, 
\qquad \omega_3(\x)=  \int_{\M} \overset{=}{R}_{\delta}(\x,\y) d \mu_{\x}.
\end{equation}
Then we have $p_{\delta}  \in H^1(\M)$, with the estimate
\begin{equation}
\left \|    p_{\delta}  \right \|_{H^1(\M)}  \leq C  ( \left \| r_{1\delta} \right \|_{H^k(\M)}+ \delta^2   \left \| \nabla_{\M} r_{\delta}  \right \|_{L^2(\M)}).
\end{equation}
\end{lmm}
This lemma guarantees the $H^1$-regularity of our nonlocal solution. For completeness, we provide its proof in Appendix \ref{appena}.

Applying Lemma \ref{higher1} directly to the bilinear form \eqref{bilinear1} yields:
\begin{equation} \label{higher20}
\left \|    u_{\delta}  \right \|_{H^1(\M)}  \leq C  ( \left \| f_{\delta} \right \|_{L^2(\M)}+ \delta^2   \left \| \nabla_{\M} f_{\delta}  \right \|_{L^2(\M)}).
\end{equation}
While the $L^2$ norm of $f_{\delta}$ was bounded in \eqref{higher15}, we now analyze  $\nabla_{\M} f_{\delta}$  through the following computation:
\begin{equation}
\begin{split}
& |\nabla_{\M}  f_{\delta}(\x)|=  \big| \nabla_{\M}^{\x} \int_{\M} f(\y) \bar{R}_{\delta}(\x,\y) d \mu_\y - \nabla_{\M}^{\x} \int_{\partial \M} (\x-\y) \cdot  \n(\y)   f(\y)      \bar{R}_{\delta} (\x, \y) d \tau_\y \big| \\
= &  \frac{1}{2\delta^2} \big| \int_{\M} f(\y)  \mathcal{T}_{\M}^{\x}(\x-\y) {R}_{\delta}(\x,\y) d \mu_\y
- 2\delta^2 \int_{\partial \M} \mathcal{T}_{\M}^{\x}  \n(\y)   f(\y)      \bar{R}_{\delta} (\x, \y) d \tau_\y 
 -  \int_{\partial \M} (\x-\y) \cdot  \n(\y)   f(\y)   \mathcal{T}_{\M}^{\x}(\x-\y)    {R}_{\delta} (\x, \y) d \tau_\y \big| \\
\leq & C ( \delta^{-1} \int_{\M} |f(\y)| {R}_{\delta}(\x,\y) d \mu_\y + \int_{\partial \M} |f(\y)| \bar{R}_{\delta}(\x,\y) d \tau_\y
+ \int_{\partial \M} |f(\y)| {R}_{\delta}(\x,\y) d \tau_\y),
\end{split}
\end{equation}
where $\mathcal{T}_{\M}^{\x}: \mathbb{R}^d \to \mathbb{R}^d$ denotes the orthogonal projection of a vector field onto the tangent space of $\M$ at $\x$. Employing the same method developed for \eqref{higher15}, we derive
\begin{equation} \label{higher19}
\begin{split}
\left \|   \nabla_{\M} f_{\delta}  \right \|_{L^2(\M)} \leq C ( \delta^{-1} \left \|    f  \right \|_{L^2(\M)} +\delta^{-1/2}  \left \|    f  \right \|_{L^2(\partial \M)} )\leq C  \delta^{-1} \left \|   f  \right \|_{H^1(\M)}.
\end{split}
\end{equation}
Substituting the bounds from \eqref{higher15} and \eqref{higher19} into \eqref{higher20}, we obtain the desired estimate:
\begin{equation}
\left \|    u_{\delta}  \right \|_{H^1(\M)}  \leq C \left \|   f  \right \|_{H^1(\M)}.
\end{equation}

\section{Vanishing Nonlocality}  \label{r2control}
This section is devoted to the proof of Theorem \ref{theorem1}(ii). Setting
\begin{equation} \label{edelta}
e_{\delta}(\x)=u(\x)-u_{\delta}(\x), \ \x \in \M; \qquad \hat{e}_{\delta}(\y)=\hat{u}(\y)-\hat{u}_{\delta}(\y)=\frac{1}{\hat{\omega}_{\delta}(\y)}
 \int_{\M} (u(\x)- u_{\delta}(\x)) \zeta_{\delta}(\x,\y) d \mu_{\x}, \ \y \in \partial \M.
\end{equation}
We observe that $e_{\delta} \in L_0^2(\M) \cap H^1(\M)$. Subtracting the first equation of \eqref{c01} from \eqref{b0019} yields
\begin{equation}
\begin{split}
 \frac{1}{\delta^2} \int_\M (e_{\delta}(\x)-e_{\delta}(\y))  R_{\delta}(\x,\y)  d \mu_\y  +   \int_{\partial \M} 
   \Delta_{ \partial \M} \hat{u}_{\delta}(\y)   \zeta_{\delta}(\x,\y)  \frac{\delta C_R}{\hat{\omega}_{\delta}(\y)}   + \Delta_{ \partial \M} u(\y)   (\x-\y) \cdot  \n(\y)  \bar{R}_{\delta} (\x, \y)   d \tau_\y 
=r_1(\x)+\tilde{f}_{\delta},
\end{split}
\end{equation}
where $r_1$ and $\tilde{f}_{\delta}$ are defined in \eqref{b0017} and \eqref{preli8} respectively. 
Now let us define two additional error functions:
\begin{equation} \label{rr2}
r_2(\x)= -\int_{\partial \M} 
   \Delta_{ \partial \M} u(\y) \big( \zeta_{\delta}(\x,\y)  \frac{\delta C_R}{\hat{\omega}_{\delta}(\y)} +   (\x-\y) \cdot  \n(\y)  \bar{R}_{\delta} (\x, \y)   \big) d \tau_\y ;
\end{equation}
\begin{equation} \label{rr3}
r_3(\x)= \int_{\partial \M} 
   \Delta_{ \partial \M} (u(\y)-\hat{u}(\y) )  \zeta_{\delta}(\x,\y)  \frac{\delta C_R}{\hat{\omega}_{\delta}(\y)}  d \tau_\y.
\end{equation}
Combining \eqref{edelta}$-$\eqref{rr3}, $e_{\delta}$ satisfies
 \begin{equation}  \label{e01}
 \begin{cases}
\displaystyle  \frac{1}{\delta^2}  \int_\M (e_{\delta}(\x)-e_{\delta}(\y)) R_{\delta}(\x,\y)  d \mu_\y  -  \delta C_R \int_{\partial \M} \Delta_{ \partial \M} \hat{e}_{\delta}(\y)   \zeta_{\delta}(\x,\y)  \frac{1}{\hat{\omega}_{\delta}(\y)} d \tau_\y 
= r_1(\x)+r_2(\x)+r_3(\x)+\tilde{f}_{\delta},  \qquad & \x \in \M; \\
 \displaystyle  \hat{e}_{\delta}(\y)- \frac{1}{\hat{\omega}_{\delta}(\y)} \int_{\M} e_{\delta}(\x) \zeta_{\delta}(\x,\y) d \mu_{\x}=0 , \qquad & \y \in \partial \M; \\
\displaystyle \int_{\M} e_{\delta}(\x) d \mu_{\x}=0,
 \end{cases}
\end{equation}
hence 
\begin{equation} \label{vanish3}
B[ e_{\delta}, q_{\delta}] =(r_1+r_2+r_3+\tilde{f}_{\delta},q_{\delta}), \qquad \forall \ q_{\delta} \in L_0^2(\M).
\end{equation}
According to Lemma \ref{higher1}, when \eqref{vanish3} holds, establishing 
\begin{equation} \label{vanish1}
|\tilde{f}_{\delta}| \leq C\delta^2  \left \| u \right \|_{H^5(\M)}, \qquad
\left \| \nabla_{\M} r_k  \right \|_{L^2(\M)}  \leq C \delta^{\frac{1}{2}}\left \| u \right \|_{H^5(\M)}, \qquad k=1,2,3;
\end{equation}
and
\begin{equation} \label{vanish2}
 (r_k, q_\delta) \leq C \delta^2 \left \| u \right \|_{H^5(\M)} (\left \| q_{\delta} \right \|_{H^1(\M)} + \left \| \bar{q}_{\delta} \right \|_{H^1(\M)}  + \left \| \overset{=}{q}_{\delta} \right \|_{H^1(\M)} ),
\ \forall \ q_{\delta} \in H^1(\M) \qquad k=1,2.3.
\end{equation} will yield
\begin{equation} \label{vanishres}
\left \| e_{\delta} \right \|_{H^1(\M)} \leq C \delta^2 \left \| u \right \|_{H^5(\M)} + C \delta^{2} \sum \limits_{r=1}^3 \left \| \nabla_{\M}r_k \right \|_{L^2(\M)} \leq  C \delta^2 \left \| u \right \|_{H^5(\M)},
\end{equation}
which precisely the part $(ii)$ of Theorem \ref{theorem1}. Thus, our remaining work is to establish \eqref{vanish1} and \eqref{vanish2}.

Since $\frac{\partial u}{\partial \n}=0$ in our setting, the term $r_1$ defined in \eqref{b0017} coincides exactly with $r_{in}$ from \eqref{b0016}. The analysis of $r_{in}$ was the focus of \cite{Base4}, where Theorem 3.1 establishes that it satisfies both \eqref{vanish1} and \eqref{vanish2}. Consequently, we need only consider the remaining terms $r_2, r_3$ and $\tilde{f}_{\delta}$.

We begin by proving that $r_2$ verifies estimates \eqref{vanish1} and \eqref{vanish2}. 
For $(\x,\y) \in \M \times \partial \M$, let us first calculate
\begin{equation}  \label{e02}
\begin{split}
& \zeta_{\delta}  (\x,\y)   \frac{\delta C_R}{\hat{\omega}_{\delta}(\y)} + (\x-\y) \cdot  \n(\y)  \bar{R}_{\delta} (\x, \y)
 = \zeta_{\delta}(\x,\y) ( \frac{\delta C_R}{\hat{\omega}_{\delta}(\y)}  -1 )
+ \big( \zeta_{\delta}(\x,\y)+ (\x-\y) \cdot  \n(\y)  \bar{R}_{\delta} (\x, \y)  \big) \\
& =\big[  \n(\y) \cdot \mathcal{T}_{\M}^{\x}(\x-\y) ( \frac{\delta C_R}{\hat{\omega}_{\delta}(\y)}  -1)  
+  ( \n(\y)- \mathcal{T}_{\M}^{\x} \n(\y)  ) \cdot ( (\x-\y)-\mathcal{T}_{\M}^{\x}(\x-\y)  ) \big]
 \bar{R}_{\delta} (\x, \y)  
=  \beta_{\delta}(\x,\y)   \bar{R}_{\delta} (\x, \y)   .
\end{split}
\end{equation} 
Then \eqref{rr2} is simplified to $r_2(\x)=-\int_{\partial \M}  \Delta_{ \partial \M} u(\y) \beta_{\delta}(\x,\y)  \bar{R}_{\delta} (\x, \y)   d \tau_\y$.
We then analyze $\hat{\omega}_{\delta}$. From its definition in \eqref{defome1}, we have
\begin{equation}  \label{e10}
\begin{split}
|\hat{\omega}_{\delta}(\y) & -\delta C_R|
=|\int_{\M} \zeta_{\delta}(\x,\y) \, d\mu_{\x} -\delta C_R |
\leq  | 2\delta^2 \int_{\partial \M}  \overset{=}{R}_{\delta}(\x,\y)  \n(\y) \cdot \n(\x) d\tau_{\x} -\delta C_R | 
 \\
& \leq \delta \ | 2\delta \int_{\partial \M} \overset{=}{R}_{\delta}(\x,\y)  d\tau_{\x} -C_R| 
+ 2\delta^2 \int_{\partial \M} (1-\n(\y) \cdot \n(\x)) \overset{=}{R}_{\delta}(\x,\y)  d\tau_{\x} 
\leq C\delta^3, \qquad \forall \ \y \in \partial \M;
\end{split}
\end{equation}
where the bound $| 2\delta \int_{\partial \M} \overset{=}{R}_{\delta}(\x,\y)  d\tau_{\x} -C_R| \leq C\delta^2$ is given in Lemma \ref{lmm3} in Appendix \ref{appenb}. The formulation \eqref{e02} thereby yields the estimate
\begin{equation} 
|\beta_{\delta}(\x,\y)| \leq C\delta^3, \qquad | \nabla_{\M}^{\x} \beta_{\delta}(\x,\y) | \leq C \delta^2,
\end{equation}
and hence
\begin{equation} \label{e06}
\begin{split}
 | \nabla_{\M}^{\x} \big( \beta_{\delta}(\x,\y) \bar{R}_{\delta} (\x, \y) \big) |
\leq C\delta^2 (   {R}_{\delta} (\x, \y) +   \bar{R}_{\delta} (\x, \y)  ) \qquad \x \in \M, \ \y \in \partial \M.
\end{split}
\end{equation} 
Therefore, using the same argument as the control of the second term in \eqref{higher15}, we obtain
\begin{equation} 
\begin{split}
& \left \| r_2 \right \|_{L^2(\M)} 
 \leq C \delta^3 \left \| \int_{\partial \M} | \Delta_{\partial \M} u(\y) | \bar{R}_{\delta} (\x, \y) d \tau_{\y}  \right \|_{L_{\x}^2(\M)} 
 \leq C\delta^{\frac{5}{2}} \left \| \Delta_{\partial \M} u \right \|_{L^2(\partial \M)} 
 \leq C\delta^{2} \left \| u \right \|_{H^3(\M)},
 \\
& \left \| \nabla_{\M} r_2 \right \|_{L^2(\M)} 
 =  \left \| \int_{\partial \M} \Delta_{\partial \M} u(\y) \big(\nabla_{\M}^{\x} \beta_{\delta}(\x,\y) \bar{R}_{\delta} (\x, \y) \big) d \tau_{\y}  \right \|_{L_{\x}^2(\M)} 
 \leq C\delta^{\frac{3}{2}} \left \| \Delta_{\partial \M} u \right \|_{L^2(\partial \M)} 
 \leq C\delta^{\frac{1}{2}} \left \| u \right \|_{H^3(\M)};
\end{split}
\end{equation} 
and hence
\begin{equation}  
\begin{split}
(r_2,q_{\delta}) \leq  C \left \| q_{\delta} \right \|_{L^2(\M)} \left \|  r_2 \right \|_{L^2(\M)} \leq C  \delta^2 \left \| q_{\delta} \right \|_{H^1(\M)} \left \|  u \right \|_{H^3(\M)}
\quad \forall \ q_{\delta} \in H^1(\M),
\end{split}
\end{equation} 
Now the proof of bound \eqref{vanish1} \eqref{vanish2} for $r_2$ is complete.  

Next, we concentrate on the control of $r_3$ defined in \eqref{rr3}.
Applying the estimate 
$$  | \nabla_{\M}^{\x} \zeta_{\delta}(\x,\y) | \leq C   (   {R}_{\delta} (\x, \y) +   \bar{R}_{\delta} (\x, \y) ) \qquad \x \in \M, \ \y \in \partial \M $$
and utilizing again the same method as the control on the second term of \eqref{higher15}, we deduce
\begin{equation}  \label{e07}
\begin{split}
 \left \| \nabla_{\M} r_3 \right \|_{L^2(\M)} 
 =  \left \| \int_{\partial \M} \Delta_{\partial \M} (u(\y)-\hat{u}(\y))  \frac{\delta C_R}{\hat{\omega}_{\delta}(\y)} (\nabla_{\M}^{\x} \zeta_{\delta}(\x,\y)) d \tau_{\y}  \right \|_{L_{\x}^2(\M)} 
 \leq C\delta^{-\frac{1}{2}} \left \|  \Delta_{\partial \M} (u-\hat{u}) \right \|_{L^2(\partial \M)} .
\end{split}
\end{equation} 
In addition, we apply integration by parts and using the same argument as \eqref{higher15} again to find
\begin{equation}  \label{e08}
\begin{split}
 (r_3& ,q_{\delta})
= \int_{\partial \M} \frac{\delta C_R}{\hat{\omega}_{\delta}(\y)}  \Delta_{ \partial \M} (u(\y) -\hat{u}(\y)) 
\int_{\M} 2\delta^2 q_{\delta}(\x) \nabla^{\x}_{\M}  \overset{=}{R}_{\delta}(\x,\y) \cdot \n(\y) d \mu_{\x} d \tau_\y \\
 = & \int_{\partial \M} \frac{\delta C_R}{\hat{\omega}_{\delta}(\y)}  \Delta_{ \partial \M} (u(\y) -\hat{u}(\y)) 
\Big( \int_{\partial \M} 2\delta^2 q_{\delta}(\x)  \n(\x) \cdot \n(\y) \overset{=}{R}_{\delta}(\x,\y) d \tau_{\x}
- \int_{\M} 2\delta^2 \nabla_{\M} q_{\delta}(\x)  \n(\y) \overset{=}{R}_{\delta}(\x,\y)  d \mu_{\x} \Big) d \tau_\y \\
 = & 2\delta^2  \int_{\partial \M} \frac{\delta C_R}{\hat{\omega}_{\delta}(\y)}  \Delta_{ \partial \M} (u(\y) -\hat{u}(\y)) \overset{=}{q}_{\delta}(\y) d \tau_{\y}
-\int_{\M} \nabla_{\M} q_{\delta}(\x) \cdot  \int_{\partial \M} \frac{2\delta^3 C_R \n(\y)}{\hat{\omega}_{\delta}(\y)}  \Delta_{ \partial \M} (u(\y) -\hat{u}(\y))  \overset{=}{R}_{\delta}(\x,\y) d \tau_{\y} d \mu_{\x} \\
& -  \int_{\partial \M} q_{\delta}(\x)  \int_{\partial \M} \frac{ 2\delta^3  C_R}{\hat{\omega}_{\delta}(\y)}  \Delta_{ \partial \M} (u(\y) -\hat{u}(\y)) 
( 1-\n(\x) \cdot \n(\y)) \overset{=}{R}_{\delta}(\x,\y) d \tau_{\y} d \tau_{\x} \\
 \leq & C \delta^2 \left \|  \Delta_{\partial \M} (u-\hat{u}) \right \|_{L^2(\partial \M)} \left \|  \overset{=}{q}_{\delta} \right \|_{L^2(\partial \M)} + C\delta^2 \left \| \nabla_{\M} q_{\delta} \right \|_{L^2(\M)} 
\left \| \int_{\partial \M} | \Delta_{ \partial \M} (u(\y) -\hat{u}(\y))|  \overset{=}{R}_{\delta}(\x,\y) d \tau_{\y} \right \|_{L^2_{\x}(\M)} \\
& +C\delta^4 \left \|  q_{\delta} \right \|_{L^2(\partial \M)} \left \| \int_{\partial \M} | \Delta_{ \partial \M} (u(\y) -\hat{u}(\y))|  \overset{=}{R}_{\delta}(\x,\y) d \tau_{\y} \right \|_{L^2_{\x}(\partial \M)} \\
 \leq & C\delta^{\frac{3}{2}} (\left \| {q}_{\delta} \right \|_{H^1(\M)} +   \left \|  \overset{=}{q}_{\delta} \right \|_{H^1(\M)}) 
 \left \|  \Delta_{\partial \M} (u-\hat{u}) \right \|_{L^2(\partial \M)}, \qquad \forall \ q_{\delta} \in H^1(\M).
   \end{split}
\end{equation}
Therefore, according to \eqref{e07} and \eqref{e08}, to prove that $r_3$ meet the conditions \eqref{vanish1} and \eqref{vanish2}, it is suffice to show
\begin{equation} \label{e09}
 \left \|  \Delta_{\partial \M} (u-\hat{u}) \right \|_{L^2(\partial \M)} \leq C \delta \left \|  u \right \|_{H^5(\M)}.
\end{equation}
Recall the definition of $\hat{u}$ in \eqref{edelta}:
\begin{equation}
 u-\hat{u}=- \frac{1}{\hat{\omega}_{\delta}(\y)} \int_{\M} (u(\x)- u(\y)) \zeta_{\delta}(\x,\y) d \mu_{\x}.
\end{equation}
We denote $s_{\delta}(\y)=\int_{\M} (u(\x)- u(\y)) \zeta_{\delta}(\x,\y) d \mu_{\x}$, and calculate
\begin{equation}
\begin{split}
 \Delta_{\partial \M} (u-\hat{u}) =& -\Delta_{\partial \M} \frac{s_{\delta}(\y)}{ \hat{\omega}_{\delta} (\y)}=-\mbox{div}_{\partial \M} (\nabla_{\partial \M} \frac{s_{\delta}(\y)}{ \hat{\omega}_{\delta} (\y)} ) =\mbox{div}_{\partial \M} ( \frac{\nabla_{\partial \M}  s_{\delta}(\y)}{ \hat{\omega}_{\delta} (\y)} -  \frac{ s_{\delta}(\y) \nabla_{\partial \M}  \hat{\omega}_{\delta}(\y) }{ \hat{\omega}^2_{\delta} (\y)}  ) \\
=& -\frac{\Delta_{\partial \M} s_{\delta}(\y)}{ \hat{\omega}_{\delta} (\y)}  +  \frac{2 \nabla_{\partial \M}s_{\delta}(\y) \cdot \nabla_{\partial \M}  \hat{\omega}_{\delta}(\y) }{ \hat{\omega}^2_{\delta} (\y)}   -s_{\delta}(\y) \frac{ \nabla_{\partial \M} \omega_{\delta} (\y) \cdot \nabla_{\partial \M} \omega_{\delta} (\y) -  \omega_{\delta} (\y)  \Delta_{\partial \M} \omega_{\delta} (\y) }{\omega^3_{\delta} (\y)}.
\end{split}
\end{equation}
Therefore, in addition to the estimate \eqref{e10}, proving the following bounds will yield \eqref{e09}:
\begin{equation} \label{e090}
\begin{split}
& |  \nabla_{\partial \M} \omega_{\delta} (\y) |  +| \Delta_{\partial \M} \omega_{\delta} (\y)| \leq C \delta^2, \quad \forall \ \y \in \partial \M;  \\
& \delta  \left \|  s_{\delta}  \right \|_{L^2(\partial \M)} + \delta \left \|   \nabla_{\partial \M} s_{\delta}  \right \|_{L^2(\partial \M)}  + \left \| \Delta_{\partial \M} s_{\delta} \right \|_{L^2(\partial \M)} \leq C\delta^2 \left \|  u \right \|_{H^5(\M)}.
\end{split}
\end{equation}
For $\omega_{\delta}$, we apply the decomposition \eqref{b0016} on the function $\n: \partial \M \to \mathbb{R}^d$ to write
\begin{equation}   \label{e091}
\begin{split}
\hat{\omega}_{\delta}(\y) =& 2\delta^2 \int_{\M}\nabla^{\x}_{\M}  \overset{=}{R}_{\delta}(\x,\y) \cdot \n(\y) d \mu_{\x}
= 2\delta^2 \int_{\partial \M} \overset{=}{R}_{\delta}(\x,\y)  d\tau_{\x} 
- 2\delta^2 \n(\y) \cdot \int_{\partial \M} (\n(\y)-\n(\x)) \overset{=}{R}_{\delta}(\x,\y)  d\tau_{\x} \\
=& 2\delta^2 \int_{\partial \M} \overset{=}{R}_{\delta}(\x,\y)  d\tau_{\x} 
+2\delta^4 \n(\y) \cdot \Delta_{\partial \M} \n(\y) \int_{\partial \M}  \overset{\equiv}{R}_{\delta}(\x,\y)  d\tau_{\x} 
+\mathcal{O}(\delta^4),
\end{split}
\end{equation}
where the kernel 
$\overset{\equiv}{R}_{\delta}(\x, \y) =C_{\delta} \overset{\equiv}{R} \big( \frac{| \x-\y| ^2}{4 \delta^2} \big)$ with $\overset{\equiv}{R}(r)=\int_r^{+\infty} \overset{=}{R}(s) ds$. By applying Lemma \ref{lmm3} in Appendix \ref{appenb}, we are able to obtain the bound \eqref{e090} for $\hat{\omega}_{\delta}$.

The control for $s_{\delta}$ is more involved. We apply the decomposition \eqref{base} two times for $u \in H^4(\partial \M)$ and $\nabla_{\M} u \in H^4(\M)$ to obtain
\begin{equation}  \label{sdelta}
\begin{split}
 s_{\delta}  (\y)  = & 2 \delta^2  \int_{\M}  (u(\x) -u(\y))  \nabla^{\x}_{\M}  \overset{=}{R}_{\delta}(\x,\y) \cdot \n(\y)    d \mu_\x  \\
=  & 2 \delta^2  \int_{\partial \M}  (u(\x) -u(\y))  \n(\y)     \overset{=}{R}_{\delta} (\x, \y)   \cdot \n(\x)  d \tau_\x 
 -2 \delta^2  \int_{\M}  \nabla_{\M} u(\x) \cdot  \n(\y)     \overset{=}{R}_{\delta} (\x, \y)   d \mu_\x  \\
= &  2 \delta^2  \int_{\partial \M}  (u(\x) -u(\y)) \overset{=}{R}_{\delta} (\x, \y)   d \tau_\x - 2\delta^2   \int_{\M}  (\nabla_{\M} u(\x)- \nabla_{\M}  u(\y))   \overset{=}{R}_{\delta} (\x, \y)   d \mu_\x \cdot  \n(\y) \\
 & - 2\delta^2 \frac{\partial u}{\partial \n} (\y) \int_{\M}    \overset{=}{R}_{\delta} (\x, \y)   d \mu_\x      
 - 2 \delta^2  \int_{\partial \M}   (u(\x) -u(\y)) (1-\n(\x)  \cdot \n(\y)) \overset{=}{R}_{\delta} (\x, \y)   d \tau_\x  \\
= & -2 \delta^4  \int_{\partial \M}  \Delta_{\partial \M} u(\x) \overset{\equiv}{R}_{\delta} (\x, \y)   d \tau_\x 
+2\delta^4 \n(\y) \cdot (  \int_{\M}  \Delta_{\M}  (\nabla_{\M} u(\x))   \overset{\equiv}{R}_{\delta} (\x, \y)   d \mu_\x -
 2 \int_{\partial \M} \overset{\equiv}{R}_{\delta}(\x,\y) \frac{\partial (\nabla_{\M} u) }{\partial \n} (\x) d \tau_\x ) \\
 & - 2\delta^2 \frac{\partial u}{\partial \n} (\y) \int_{\M}    \overset{=}{R}_{\delta} (\x, \y)   d \mu_\x     
 - 2 \delta^2  \int_{\partial \M}   (u(\x) -u(\y)) (1-\n(\x)  \cdot \n(\y)) \overset{=}{R}_{\delta} (\x, \y)   d \tau_\x +2\delta^4 r_s(\y).
\end{split}
\end{equation}
Here $r_s(\y)$ denotes the residue terms of the decomposition \eqref{base}, which has been completely controlled in \cite{Base1}. For the remaining terms in \eqref{sdelta}, By applying Lemma \ref{lmm4}  \ref{lmm5}  \ref{lmm4}  
\ref{lmm2} and Corollary \ref{lmm6} in Appendix \ref{appenb} respectively, we are able to conclude
\begin{equation} \label{e092}
\begin{split}
&    \left \|  s_{\delta}  \right \|_{L^2(\partial \M)} +  \left \|   \nabla_{\partial \M} s_{\delta}  \right \|_{L^2(\partial \M)}  + \left \| \Delta_{\partial \M} s_{\delta} \right \|_{L^2(\partial \M)}
 \\
  \leq & C ( \delta^{\frac{5}{2}} \left \|  \Delta_{\partial \M} u \right \|_{H^2(\partial \M)} 
  + \delta^{\frac{7}{2}} \left \|  \Delta_{\M}  (\nabla_{\M} u ) \right \|_{H^2(\M)} 
  +  \delta^{\frac{5}{2}} \left \| \frac{ (\nabla_{\M} u) }{\partial \n} \right \|_{H^2(\partial \M)} 
  + \delta^2  \left \|  \frac{\partial u}{\partial \n} \right \|_{H^2(\M)} +  \delta^{\frac{5}{2}} \left \| u \right \|_{H^2(\partial \M)} ) \\
  \leq & C \delta^2  \left \|  u \right \|_{H^5(\M)}.
\end{split}
\end{equation}
This completes the second part of \eqref{e090}, which yields \eqref{e09}, and consequently $r_3$ meets the condition \eqref{vanish1} and \eqref{vanish2}.

For convenience, the control \eqref{vanish1} for $\tilde{f}_{\delta}$ have been moved to Appendix \ref{appenf}. 

We thereby obtain \eqref{vanishres} since all the truncation errors have been controlled, while \eqref{vanishres} is exactly the second part of Theorem \ref{theorem1}.

\section{Model Generalizations}
The preceding analysis focused on establishing convergence rates between the nonlocal model \eqref{c01} and its local counterpart \eqref{bg01}. We now extend this framework to three related local problems, each representing a modification of \eqref{bg01}, and develop their corresponding nonlocal formulations.

\subsection{Manifold Poisson Model}
We begin by considering the following modified problem:
\begin{equation}   \label{Poisson}
\begin{cases}
-\Delta_{\M} u(\x) +\lambda(\x) u(\x)=f(\x) & \x \in \M, \\
\frac{\partial u}{\partial \n}(\x)=0  & \x \in \partial \M,  
\end{cases}
\end{equation}
here $f$ is assumed to be $H^3(\M)$, $\lambda$ is a sufficiently smooth function on $\M$ with the assumption $\lambda>\lambda_0$ for some constant $\lambda_0>0$. This problem is \eqref{bg01} added by the $\lambda u$ term. From classical elliptic theory mentioned in \cite{Base3}, it follows that \eqref{Poisson} has a unique solution $u \in H^5(\M)$.

To derive a nonlocal approximation of \eqref{Poisson}, the most direct approach is to plug the equation $\Delta_{\M} u=-f+u$ into \eqref{b0016} instead of $\Delta_{\M}u=-f$, which yields the appearance of two additional terms involving $\lambda$:
\begin{equation}  \label{pc01}
 \begin{split}
& \frac{1}{\delta^2} \int_\M (u(\x)-u(\y))  R_{\delta}(\x,\y)  d \mu_\y   + \int_{\M}  \lambda(\y) u(\y) \bar{R}_{\delta}(\x,\y) d \mu_\y + \int_{\partial \M} -(\x-\y) \cdot \n(\y)  \lambda(\y) u(\y) \bar{R}_{\delta}(\x,\y) d \tau_{\y} \\
&+   \int_{\partial \M} \Delta_{ \partial \M} u(\y)   (\x-\y) \cdot  \n(\y)  \bar{R}_{\delta} (\x, \y)      d \tau_\y 
= f_{\delta}(\x) +r_1(\x), \ \x \in \M.
 \end{split}
\end{equation}
To ensure energy non-negativity for the additional two terms as well, we define the auxiliary variables
 \begin{equation} \label{omega0}
\bar{u}_{\delta} (\y) =\frac{ 1 }{ \omega_2(\y) } \int_{\M} u_{\delta} (\z)   \bar{R}_{\delta}(\z,\y) d \mu_{\z}, 
\qquad \omega_2(\y)= \int_{\M}   \bar{R}_{\delta}(\z,\y) d \mu_{\z} \qquad \y \in \M;
\end{equation}
 \begin{equation}
\hat{u}_{\delta} (\y) =\frac{ 1 }{ \hat{\omega}_{\delta}(\y) } \int_{\M} u_{\delta}(\z)   \zeta_{\delta}(\z,\y) d \mu_{\z}, 
\qquad \hat{\omega}_{\delta}(\y)= \int_{\M}   \zeta_{\delta}(\z,\y) d \mu_{\z} \qquad \y \in \partial \M.
\end{equation}
By mimicking the process of construction \eqref{b0017}-\eqref{preli8}, we write down the final nonlocal model that approximates \eqref{Poisson}:
 \begin{equation}  \label{pc02}
 \begin{split}
& \frac{1}{\delta^2} \int_\M (u_{\delta}(\x)-u_{\delta}(\y)) R_{\delta}(\x,\y)  d \mu_\y  +  \int_{\M} \lambda(\y)  \bar{u}_{\delta}(\y) \bar{R}_{\delta}(\x,\y) d \mu_\y 
+ \int_{\partial \M}   \lambda(\y) \hat{u}_{\delta}(\y) \zeta_{\delta}(\x,\y) d \tau_{\y} 
 \\
& -  \delta C_R \int_{\partial \M} \Delta_{ \partial \M} \hat{u}_{\delta}(\y)   \zeta_{\delta}(\x,\y)  \frac{1}{\hat{\omega}_{\delta}(\y)} d \tau_\y = f_{\delta}(\x)  \qquad \x \in \M. 
 \end{split}
\end{equation}
In such case, the mean-zero constraint is no longer imposed since the left-side operator in \eqref{pc02} is always coercive in $L^2(\M)$. Using the same idea as in Section \ref{psd}, we are able to prove similar well-posedness argument for \eqref{pc02}. To prove its second-order convergence to local problem \eqref{Poisson}, in addition to the truncation errors mentioned in Section \ref{r2control}, two other errors need to be bounded are
\begin{equation} \label{higher174}
\begin{split}
\left \| u-\bar{u} \right \|_{L^2(\M)} \leq C \delta^2  \left \| u \right \|_{H^5(\M)} , \qquad
\left \| u-\hat{u} \right \|_{L^2(\partial \M)}  \leq C \delta  \left \| u \right \|_{H^5(\M)} .
\end{split}
\end{equation}
As $u-\hat{u}=s_{\delta}/ \hat{\omega}_{\delta}$, we can directly apply \eqref{e091}$-$\eqref{e092} to give its $L^2$ bound on $\partial \M$. For the remaining term $u-\hat{u}$, we can apply the decomposition \eqref{b0016} on $u \in H^4(\M)$ to derive such bound.
Consequently, a full error analysis demonstrates that the solution $u_\delta$ of the nonlocal model \eqref{pc02} converges to the solution $u$ of the local problem \eqref{Poisson} with second-order accuracy in $\delta$.

\subsection{Model with Non-Homogeneous Boundary Condition}
We now consider the following local model subject to non-homogeneous boundary conditions:
 \begin{equation}   \label{nonhon}
\begin{cases}
-\Delta_{\M} u(\x) =f(\x) & \x \in \M, \\
\frac{\partial u}{\partial \n}(\x)=g(\x)  & \x \in \partial \M, \\
\int_{\M} u(\x) d \mu_{\x}=0.
\end{cases}
\end{equation}
where $f \in H^3(\M), g \in H^{\frac{7}{2}}(\partial \M)$, and subject to the compatibility condition $\int_{\M} f(\x) d \mu_{\x} + \int_{\partial \M} g(\x) d \tau_{\x}=0$, which ensures the well-posedness of the problem.

To derive its nonlocal approximation, we start from \eqref{b0016} again.
The equality \eqref{A.2} in Appendix \ref{appen0} gives
$$\Delta_\M u(\y) =  \Delta_{\partial \M} u(\y) +\nabla_{\partial \M}^2 u(\n,\n) (\y) - \kappa(\y) \ \frac{\partial u}{\partial \n}(\y), \qquad \forall \ \y \in \M,$$
where $\kappa(\y)= (\frac{1}{\sqrt{\det g_{\partial}}}  \frac{\partial }{\partial t} \sqrt{\det g_{\partial}}  )(\y)$ is a smooth function on $\y$ that describes the curvature of $\partial \M$ on $\M$ at $\y$, with $g_{\partial}$ and $t$ defined in Appendix \ref{appen0}. Following the same analytical procedure employed in deriving \eqref{c01}, we replace $\frac{\partial u}{\partial \n}$ by $g$ to obtain the following nonlocal approximation of problem \eqref{nonhon}:
 \begin{equation}  \label{non02}
 \begin{cases}
\displaystyle  \frac{1}{\delta^2} \int_\M (u_{\delta}(\x)-u_{\delta}(\y)) R_{\delta}(\x,\y)  d \mu_\y - \delta C_R \int_{\partial \M} \Delta_{ \partial \M} \hat{u}_{\delta}(\y)   \zeta_{\delta}(\x,\y)  \frac{1}{\hat{\omega}_{\delta}(\y)} d \tau_\y 
 = f_{\delta}(\x)+g_{\delta}(\x) - \tilde{f}_{\delta}-\tilde{g}_{\delta},   & \x \in \M; \\
 \displaystyle  \hat{u}_{\delta}(\y)- \frac{1}{\hat{\omega}_{\delta}(\y)} \int_{\M} u_{\delta}(\x) \zeta_{\delta}(\x,\y) d \mu_{\x}=0 , \qquad & \y \in \partial \M; \\
\displaystyle \int_{\M} u_{\delta}(\x) d \mu_{\x}=0,
 \end{cases}
\end{equation}

The constant functions $\tilde{f}$ and $\tilde{g}$ are introduced to ensure the right-hand side of \eqref{non02} has zero average, where $\tilde{f}_{\delta}$ is defined in \eqref{preli8}, and $\tilde{g}_{\delta}$ is defined as
\begin{equation}
\tilde{g}_{\delta}=\frac{1}{| \M |} \int_{\M}  \int_{\partial \M}  \big( 2 \bar{R}_{\delta} (\x, \y)  +  \zeta_{\delta}(\x,\y)  \big) g(\y)    d \tau_\y  d \mu_\x.
\end{equation}
It can be proved that $ \tilde{f}_{\delta}+\tilde{g}_{\delta}=\mathcal{O}(\delta^2)$ given the compatibility condition, which is a direct corollary of Appendix \ref{appenf}. Following the same reasoning as in the derivation of \eqref{c01}, the nonlocal model \eqref{non02} admits a unique solution $u_{\delta} \in H^1(\M)$ with zero mean, satisfying:
\begin{equation}
\left \| u_{\delta} \right \|_{H^1(\M)}   \leq C (\left \| f \right \|_{H^1(\M)} +  \left \| g \right \|_{H^1(\partial \M)}) ,
\end{equation}
while approximating the solution $u$ of problem \eqref{nonhon} with convergence rate
\begin{equation}
\left \| u-u_{\delta} \right \|_{H^1(\M)}   \leq C \delta^2 ( \left \| f \right \|_{H^3(\M)} + \left \| g \right \|_{H^{\frac{7}{2}}(\partial \M)}),
\end{equation}

\subsection{Model with Nonlinear Term}
In the final part of this section, we investigate the construction of a nonlocal approximation for the following local model:
 \begin{equation}  \label{nonli01}
\begin{cases}
-\Delta_{\M} u(\x) + \lambda u(\x) |u^{2p-2}(\x)|=f(\x) & \x \in \M, \\ 
\frac{\partial u}{\partial \n}(\x)=0  & \x \in \partial \M, 
\end{cases}
\end{equation}
where $f \in H^2(\M)$, $\lambda > 0$ is a constant, and $p > 1$. The presence of the nonlinear term $u|u|^{2p-2}$ introduces additional complexity to the problem. Through a straightforward variational approach, this problem can be reformulated as finding critical points of the energy functional:
 \begin{equation}
  J(u) = \frac{1}{2} \int_{\mathcal{M}} |\nabla u(\x) |^2  d \mu_{\x} + \frac{\lambda}{2p} \int_{\mathcal{M}} |u(\x)|^{2p} \, d \mu_{\x} - \int_{\mathcal{M}} f(\x) u(\x) \, d \mu_{\x}.
 \end{equation}
\begin{itemize}
\item The functional $J(u)$ is strictly coercive. When the nonlinearity is subcritical with $p < \frac{m}{m-2}$, the Sobolev embedding $H^1(\mathcal{M}) \hookrightarrow L^{2p}(\mathcal{M})$ becomes compact. This compactness ensures the weak lower semicontinuity of $J$. Consequently, $J$ admits a minimizer in $H^1(\mathcal{M})$, which corresponds to a weak solution of the problem.
\item In addition to the subcritical condition $p < \frac{m}{m-2}$, when $1<p<2$, the functional $J(u)$ becomes strictly convex, This convexity property guarantees the uniqueness of the solution to problem \eqref{nonli01}.
\item Assume $p < \frac{m}{m-2}$. The definition of $J(u)$ immediately implies that $u \in H^1(\mathcal{M}) \cap L^{2p}(\mathcal{M})$, and consequently $u|u|^{2p-2} \in L^{\frac{2p}{2p-1}}(\mathcal{M})$. Applying elliptic regularity theory to problem \eqref{nonli01}, we obtain $u \in W^{2, \frac{2p}{2p-1}}(\mathcal{M})$. If $p$ further satisfies the embedding condition $W^{2, \frac{2p}{2p-1}}(\mathcal{M}) \hookrightarrow L^r(\mathcal{M})$ for some $r > 2p$, then a bootstrap argument yields higher regularity of $u$.
\end{itemize}
Following the formulation in \eqref{c01}, we derive a nonlocal approximation of problem \eqref{nonli01} given by
 \begin{equation}  \label{nonli02}
 \begin{split}
& \frac{1}{\delta^2} \int_\M (u_{\delta}(\x)-u_{\delta}(\y)) R_{\delta}(\x,\y)  d \mu_\y  +  \int_{\M} \lambda  \bar{u}_{\delta}(\y) |  \bar{u}_{\delta}(\y)|^{2p-2} \bar{R}_{\delta}(\x,\y) d \mu_\y 
+ \int_{\partial \M}   \lambda \hat{u}_{\delta}(\y) |  \hat{u}_{\delta}(\y)|^{2p-2} \zeta_{\delta}(\x,\y) d \tau_{\y} 
\\  &
  - \delta C_R \int_{\partial \M} \Delta_{ \partial \M} \hat{u}_{\delta}(\y)   \zeta_{\delta}(\x,\y)  \frac{1}{\hat{\omega}_{\delta}(\y)} d \tau_\y 
 = f_{\delta} (\x)  \qquad \x \in \M.
 \end{split}
\end{equation}
where $\bar{u}_{\delta}, \hat{u}_{\delta}$ are defined the same as \eqref{omega0}. The corresponding energy functional is given by
\begin{equation} \label{nonli03}
\begin{split}
& J_{\delta}(u_{\delta}) =\frac{1}{\delta^2 } \int_\M \int_\M (u_{\delta}(\x)-u_{\delta}(\y))^2 R_{\delta}(\x,\y)  d \mu_\y d \mu_{\x}
+ \frac{\lambda}{2p} \int_{\M} \omega_2(\y) | \bar{u}_{\delta}(\y) |^{2p} d \mu_{\y} + \frac{\lambda}{2p} \int_{\M} \hat{\omega}_{\delta}(\y) | \hat{u}_{\delta}(\y) |^{2p} d \tau_{\y} \\
& + \delta C_R \int_{\partial \M}   \nabla_{ \partial \M} \hat{u}_{\delta}(\y)  \cdot \nabla_{ \partial \M} \hat{u}_{\delta}(\y) d \tau_{\y}
-  \int_{\M} {u}_{\delta}(\x) f_{\delta}(\x)  d\mu_{\x}
\end{split}
\end{equation}
Following the same reasoning used to bound the bilinear form $B_{\delta}[\cdot, \cdot]$
in the proofs of Lemma \ref{Lax} and Lemma \ref{higher1}, we obtain the following estimates for $J_{\delta}$:
\begin{equation}
J_{\delta}(u_{\delta}) + C_1 \left \| f \right \|^2_{H^1(\M)} \geq C_2 \left \| u_{\delta} \right \|^2_{H^1(\M)} ,
\qquad
J_{\delta}(u_{\delta}) \leq C_{\delta} (\left \| u_{\delta} \right \|^2_{L^{2p}(\M)} + \left \| f \right \|^2_{H^1(\M)} );
\end{equation}
Here, $C_{\delta} $ denotes a constant depending on $\delta $, whereas $C_1$ and $C_2$ are independent of $\delta$. The estimates above imply that $J_{\delta}$ is strictly coercive and weakly lower semicontinuous in $H^1(\mathcal{M})$, provided the exponent $p$ is subcritical, i.e., $p < \frac{m}{m-2}$.  
By applying the direct method in the calculus of variations, we conclude that $J_{\delta}$ admits a critical point, which in turn yields the existence of a solution to problem \eqref{nonli02}.

From the given conditions, we are able to deduce the bound $ \left \| u_{\delta} \right \|_{H^1(\mathcal{M})} \leq C \left \| f \right \|_{H^1(\mathcal{M})}. $
Although $ u_{\delta} \in H^k(\mathcal{M}) $ for any $ k \geq 1 $, the $ L^{2p}(\mathcal{M}) $-norm of $ u_{\delta} $ lacks a $\delta-$independent upper bound.  
Furthermore, uniqueness of the solution cannot be guaranteed, as $J_{\delta}$ is not necessarily strictly convex.  
While a rigorous vanishing-nonlocality theory for \eqref{nonli02} remains open, numerical experiments suggest a favorable convergence rate.

\section{A Numerical Example} \label{appenc}

In this part, we will conduct numerical tests to verify the convergence of our model \eqref{c01}. 
Assuming we are given the input value $\delta$, the kernel function $R(r)$, the set of points $\PP=\{ \textbf{p}_i \}_{i=1}^{n_0} $ that samples $\M$, $\QQ=\{ \textbf{q}_k\}_{k=1}^{m_0}$ that samples $\partial \M$;  the volume weight $\mathcal{A}=\{A_i \}_{i=1}^{n_0}$ for each $\textbf{p}_i  \in \M$, and the hypersurface weight $\LL=\{ L_k \}_{k=1}^{m_0}$ for each $\textbf{q}_k \in \partial \M$.

Then we discretize \eqref{c01} in the most straightforward way:
\begin{equation}  \label{numericalsoln}
\begin{cases}
\displaystyle \sum \limits_{j=1}^{n_0} \frac{1}{\delta^2} R_\delta(\textbf{p}_i, \textbf{p}_j) (u_i-u_j) A_j  - \delta C_R \sum_{k=1}^{m_0}   (\Delta^0_{\partial \M} v_k )     \zeta_{\delta}(\textbf{p}_i, \textbf{q}_k) \frac{1}{\hat{\omega}_{\delta} (\textbf{q}_k)}    L_k    = f^i_{\delta} -\tilde{f}^0_{\delta}, & i=1,2,...n_0; \\
\displaystyle   v_k- \frac{1}{\hat{\omega}_{\delta} (\textbf{q}_k)}  \sum \limits_{i=1}^{n_0} u_i \zeta_{\delta}(\textbf{p}_i, \textbf{q}_k) A_i   =0 , & k=1,2,...,m_0; \\
 \sum \limits_{i=1}^{n_0} u_i A_i=0.
\end{cases}
\end{equation}
Here $R_{\delta}$ and $\bar{R}_{\delta}$ are the kernel functions defined in Section 1, and
\begin{equation}
\begin{split}
\zeta_{\delta}(\textbf{p}_r, \textbf{q}_k)=-( \textbf{p}_r - \textbf{q}_k) \cdot \n( \textbf{q}_k ) \bar{R}_{\delta}(\textbf{p}_r,  \textbf{q}_k ), 
\qquad
 \hat{\omega}_{\delta} (\textbf{q}_k)=\sum \limits_{r=1}^{n_0}  \zeta_{\delta}(\textbf{p}_r, \textbf{q}_k) A_r,
\end{split}
\end{equation}
Besides, $ \Delta^0_{\partial \M}$ is the discrete form of $\Delta_{\partial \M}$, where we rely on \eqref{b0016} to construct it as
\begin{equation} 
\begin{split}
 \Delta^0_{\partial \M} v_k=   -   \frac{ 2 }{  \delta C_R } \sum \limits_{l=1}^{m_0} ( v_k-v_l) \bar{R}_{\delta}( \textbf{q}_k, \textbf{q}_l) L_l,
\end{split}
\end{equation}
so that we can avoid the calculation of $C_R$ in the model. The source vector is defined by
\begin{equation}
f_{\delta}^i=\sum \limits_{j=1}^{n_0} f(\textbf{p}_j) \bar{R}_\delta( \textbf{p}_i, \textbf{p}_j) A_j + \sum \limits_{k=1}^{m_0} f(\textbf{q}_k) \zeta_\delta( \textbf{p}_i, \textbf{q}_k) L_k,
\qquad
\tilde{f}^0_{\delta}=(\sum \limits_{i=1}^{n_0} f_{\delta}^i A_i ) / (\sum \limits_{i=1}^{n_0}  A_i). 
\end{equation}
By observation, setting $L_k=0$ for $k=1,...,m_0$ will transform \eqref{numericalsoln} into the reduced numerical model without the boundary terms, as proposed in eq.10 of \cite{Base1}.
Let us first prove that \eqref{numericalsoln} assures a unique solution vector $(U,V)= ( \{u_i\}_{i=1,...,n_0}, \{v_k\}_{k=1,2,...,m_0})$. The first two lines of system \eqref{numericalsoln} give $(m_0+n_0)$ linear equations on the variables $\{u_i\}_{i=1,...,n_0}$, $\{v_k\}_{k=1,2,...,m_0}$,  in the matrix form:
\begin{equation} \label{matrix1}
\begin{pmatrix}
 \frac{1}{\delta^2} A^{-1} {R_A} & 2 \bm{\zeta}  \bar{R}_L \\
- \bm{\zeta}^T A  &  I_{n_0}
\end{pmatrix}
\begin{pmatrix}
U \\
V
\end{pmatrix}
=
\begin{pmatrix}
F \\
\textbf{0}
\end{pmatrix},
\end{equation}
where $A=\mbox{diag} \{ A_1, ... ,A_{n_0} \},  \bm{\zeta}=\{ \frac{ \zeta_{\delta}(\textbf{p}_r, \textbf{q}_k)}{\hat{\omega}_{\delta}(\textbf{q}_k) }  \}_{n_0 \times m_0} $; $R_A$ and $\bar{R}_L$ are given respectively by:
\begin{equation} \label{matrix2b}
{R_A}_{ij}=
\begin{cases}
R_{\delta} (\textbf{p}_i, \textbf{p}_j) A_i A_j &   i \neq j  \\
\sum_{j' \neq i} R_{\delta} (\textbf{p}_i, \textbf{p}_{j'}) A_i A_{j'}  & i=j
\end{cases}
;
\qquad
\bar{R}_{Lkl}=
\begin{cases}
\bar{R}_{\delta} (\textbf{p}_k, \textbf{p}_l) L_k L_l &   k \neq l  \\
\sum_{l' \neq k} \bar{R}_{\delta} (\textbf{p}_k, \textbf{p}_{l'} ) L_k L_{l'}  & k=l
\end{cases}.
\end{equation}
where $1\leq i,j \leq n_0, \ 1 \leq q,k \leq m_0.$ 
Besides, the source vector $F=\{f_{\delta}^i-\tilde{f}_{\delta}^0\}_{i=1,...,n_0}$. By plugging in $V=-\bm{\zeta}^T AU $, we further write \eqref{matrix1} as
\begin{equation} \label{stiff3}
  ( \frac{1}{\delta^2} R_A + 2 A \bm{\zeta}  \bar{R}_L \bm{\zeta}^T A^T) U= AF .
 \end{equation}
The stiff matrix in \eqref{stiff3} is symmetric, semi-positive definite, with null space 
$S=\{c \textbf{1} | c \in \mathbb{R} \}$, where $\textbf{1}$ denotes the vector of all ones. The RHS term $AF$ lies in $S^{\perp}$. With the third equation of \eqref{numericalsoln} imposed, \eqref{stiff3} assures a unique solution vector $U$.

Our purpose next is to compare $U$ with the exact solution $u$ that solves the local problem \eqref{bg01}. Since the $H^1$ error is difficult to compute, here we only show the $L_2$ error
 \begin{equation} \label{ee2}
e_2=\sqrt{ \frac{
 \sum \limits_{j=1}^{n_0} ( u_{j}-u(\textbf{p}_{j}))^2 A_{j} }{ \sum \limits_{j=1}^{n_0}  u^2(\textbf{p}_j) A_{j} }},
 \end{equation}
In addition, for the numerical examples in this section, we choose the following kernel function $R$ for convenience:
\begin{equation} \label{kernelr}
R(r)=
\begin{cases}
\frac{1}{2}(1+ \cos \pi r ), & 0 \leq r \leq 1, \\
0, & r>1.
\end{cases}
\end{equation}
so that the functions $R_{\delta}, \bar{R}_{\delta}$ and $\overset{=}{R}_{\delta}$ in \eqref{numericalsoln} can be calculated through the definitions \eqref{kernelori} when $\delta$ is given. 

\subsection{Example: Unit Hemisphere $S^2$ in $\mathbb{R}^3$ } \label{hmsphere}
In the first example, we let the manifold $\M$ be the upper half of the unit hemisphere, with equation
\begin{equation}
x^2+y^2+z^2=1, \qquad z \geq \frac{1}{2},
\end{equation}
with its boundary $\partial \M$ be the unit circle $x^2+y^2=1, z =\frac{1}{2}$. 
In the local Poisson problem \eqref{bg01}, we let the exact solution $u$ be $u(x,y,z)=(z-\frac{1}{2})^2-\frac{1}{12}$. By the definition of $\nabla_{\M}$ and $\Delta_{\M}$, we parametrize $\M$ and calculate 
\begin{equation}
\begin{split}
& -\Delta_{\M} u(x,y,z)=6z^2-2z-2, \qquad   (x,y,z) \in \M,     \\
& \frac{\partial u}{\partial \n}(x,y, z) =0,  \qquad  (x,y,z) \in \partial \M, \\
& \int_{\M} u \ d \mu_{\x}=0.
 \end{split}
\end{equation}
The volume weights are calculated as follows: \\
for $t=5:5:90$
\begin{enumerate} 
\item $\delta=\sqrt{1/t}, \ n_0=t^2+3t, \ m_0=3t$;
\item generate random matrix $B=\mbox{rand}(n_0,2)$. \\
for $i \in [1, n_0-m_0]$,
\begin{equation} \label{matrix3}
\textbf{p}_i=
\begin{pmatrix}
\sqrt{1- (\frac{B(i,2)+1}{2})^2} \cos (2\pi B(i,1)), \ &
 \ \sqrt{1- (\frac{B(i,2)+1}{2})^2} \sin(2\pi B(i,1)), \ &
\ \frac{B(i,2)+1}{2}
\end{pmatrix};
\end{equation}
for $i \in [n_0-m_0+1,n_0]$
\begin{equation} \label{matrix4}
\textbf{p}_i=
\begin{pmatrix}
\frac{\sqrt{3}}{2} \cos (2\pi B(i,1)) , \ &
 \  \frac{\sqrt{3}}{2} \sin(2\pi B(i,1)), \ &
\ \frac{1}{2} 
\end{pmatrix};
\end{equation}
$\textbf{q}_k=\textbf{p}_{k+n_0-m_0}, k=1,2,...,m_0$;
\item for $i \in [1,n_0]$, find $20$ nearest neighbors of $\textbf{p}_i$ in $\PP$, project them onto the tangent plane of $\M$ at $\textbf{p}_i$; construct a $2d$ triangulation of these 21 points using delaunay triangulation;
\item $A_i$ is assigned to be $1/3$ the sum of areas of triangles with $\textbf{p}_i$ as vertex;
\item In full model, $ L_k= \frac{1}{2}( \big\vert \textbf{q}_{k_1}- \textbf{q}_{k} \big\vert+\big\vert \textbf{q}_{k}-\textbf{q}_{k_{2}} \big\vert), k \in [1,m_0] $, where $\textbf{q}_{k_1}, \textbf{q}_{k_2}$  are two nearest neighbors of $\textbf{q}_{k}$ in boundary points. In reduced model, $L_k=0, k \in [1,m_0]$.
\end{enumerate}
   
\begin{figure}[htbp] 
\centering
\includegraphics[width=\textwidth]{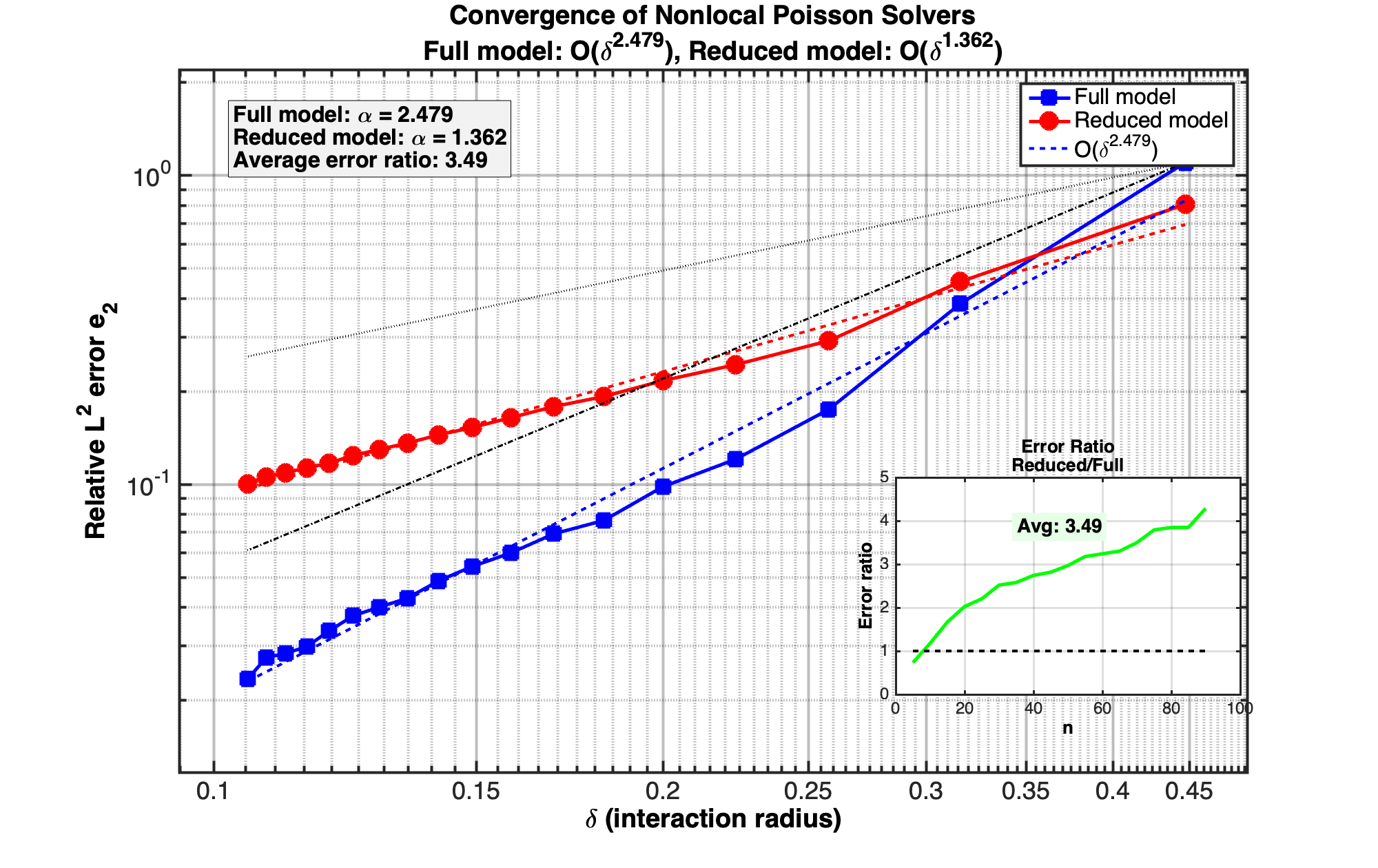}
\caption{Numerical Example \ref{hmsphere}}
\label{figure1}
\end{figure}

The results are presented in Fig. \ref{figure1}. As shown, $e_2$ exhibits an approximately linear dependence on $\delta^{2.5}$-an improvement of nearly one order of magnitude compared to the old model without boundary terms. The observed deviations are primarily attributed to the inherent randomness of the point cloud.

\subsection{Example: Unit 3-Hemisphere $S^3$ in $\mathbb{R}^4$ } \label{ball}
In the second example, let $\M$ be the manifold embedded in $\mathbb{R}^4$:
\begin{equation}
x^2+y^2+z^2+w^2=1, \qquad w \geq 0,
\end{equation}
with its boundary $\partial \M$ be the unit ball $x^2+y^2+z^2=1, w =0$. 
Similarly, let $u(x,y,z,w)=xyz$ be the exact solution of \eqref{bg01}. 
Then
\begin{equation}
\begin{split}
& - \Delta_{\M} u(x,y,z,w)=15xyz, \qquad   (x,y,z,w) \in \M,     \\
& \frac{\partial u}{\partial w}(x,y, z,w) =0,  \qquad  (x,y,z,w) \in \partial \M. \\
& \int_{\M} u \ d \mu_{\x}=0.
 \end{split}
\end{equation}

\begin{figure}[htbp] 
\centering
\includegraphics[width=\textwidth]{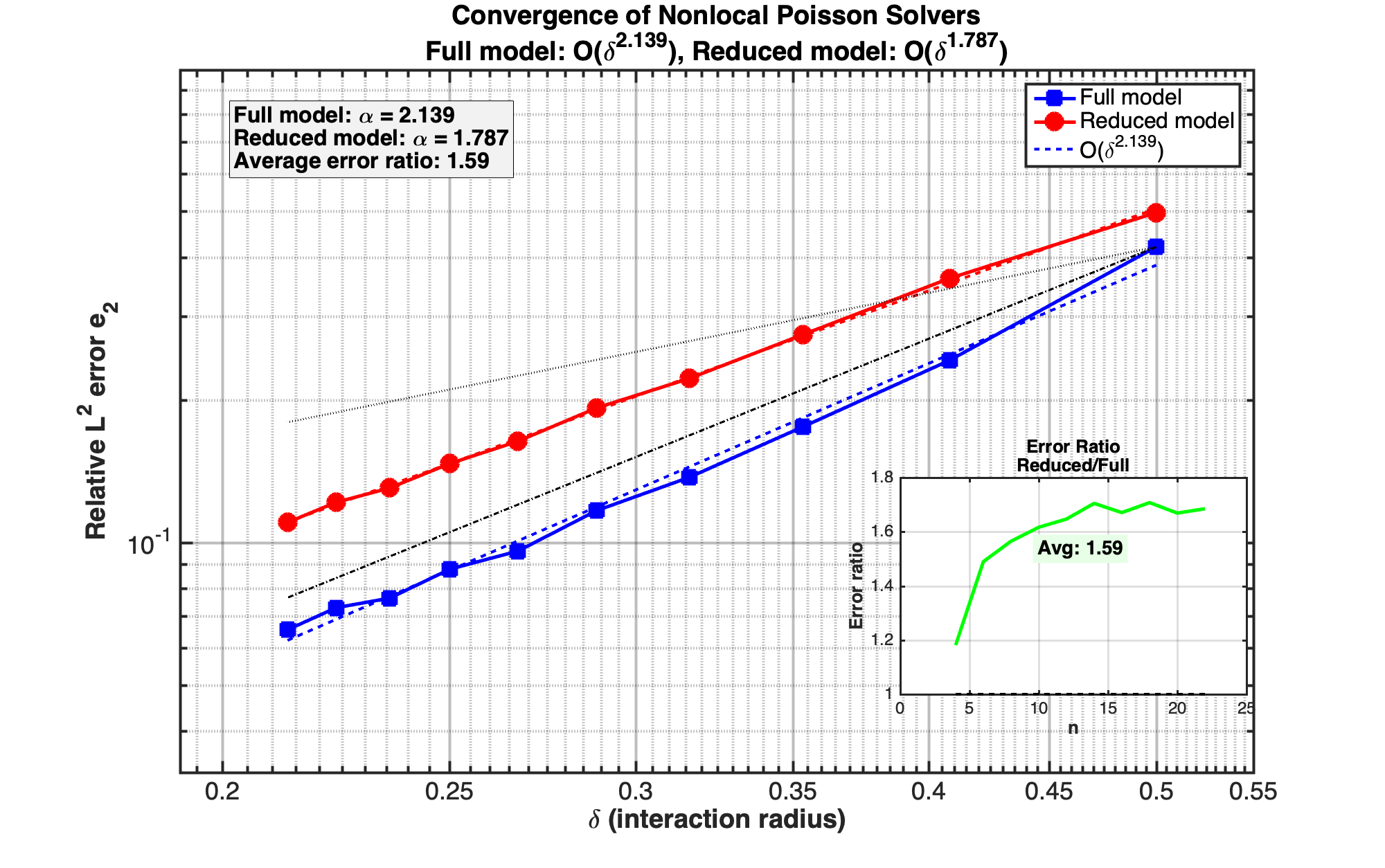}
\caption{Example \ref{ball}}
\label{figure2}
\end{figure}
The volume weights are computed as follows: \\
for $t=4:2:22$
\begin{enumerate} 
\item $\delta=\sqrt{1/t}, \ n_0=3t^3+4t^2, \ m_0=4t^2$;
\item generate $B=\mbox{rand}(n_0,3)$. \\
for $i \in [1, n_0-m_0]$,
\begin{equation} 
\textbf{p}_i=
\begin{pmatrix}
\sqrt{B(i, 1)} \cos (2\pi B(i,2)),   &
  \sqrt{B( i, 1)} \sin(2\pi B(i,2)),  &
 \sqrt{1-B(i, 1)}\cos (\pi B(i,3)), &
 \sqrt{1-B(i, 1)}\sin (\pi B(i,3)) 
\end{pmatrix};
\end{equation}
for $i \in [n_0-m_0+1,n_0]$
\begin{equation} 
\textbf{p}_i=
\begin{pmatrix}
\sqrt{1-(2B(i,2)-1)^2} \cos (2\pi B(i,1)) ,  &
  \sqrt{1-(2B(i,2)-1)^2}  \sin(2\pi B(i,1)),  &
2B(i,2)-1,   &
 0
\end{pmatrix};
\end{equation}
$\textbf{q}_k=\textbf{p}_{k+n_0-m_0}, k=1,2,...,m_0$;
\item for $i \in [1,n_0]$, find $50$ nearest neighbors of $\textbf{p}_i$ in $\PP$, project them onto the $3d$ tangent hyperplane of $\M$ at $\textbf{p}_i$; construct a tetrahedron mesh of these 51 points using Delaunay algorithm;
\item $A_i$ is assigned to be $1/4$ the sum of volumes of tetrahedrons with $\textbf{p}_i$ as vertex;
\item In full model, $L_k, k \in [1,m_0]$ is calculated in the same way as $A_i$ in the previous example. In reduced model, $L_k \equiv 0$.
\end{enumerate}

As shown in Fig. \ref{figure2}, the $3D$ case also yields convergence rates exceeding second order, indicating higher efficiency than the reduced model.

\section{Conclusion and Future Work}

In this paper, we propose a nonlocal approximation of the manifold Poisson model with Neumann boundary conditions. We analyze the well-posedness of the model and establish a second-order localization rate.
A key advantage of our nonlocal model is its compatibility with meshless numerical schemes, such as the Point Integral Method (PIM) \cite{LSS}  \cite{weightedLaplacian}. Unlike traditional Manifold Finite Element Methods (MFEM), which require a global mesh, PIM only relies on local point clouds, significantly simplifying discretization. Moreover, the numerical solution obtained from our nonlocal model can also serve as an approximation to the corresponding local Poisson problem.

In subsequent work, we plan to extend this framework to derive accurate nonlocal approximations for other PDEs, including elliptic transmission problems, Stokes equations, diffusion-reaction systems, and Poisson problems on polygonal domains.

\begin{appendix}

\section{Proof of Lemma \ref{higher1}}  \label{appena}

\begin{proof}
As discussed, \eqref{higher2} implies
\begin{equation} \label{higher4}
 \frac{1}{\delta^2}  \int_\M (p_{\delta}(\x)-p_{\delta}(\y)) R_{\delta}(\x,\y)  d \mu_\y  -  \delta C_R \int_{\partial \M} \Delta_{ \partial \M} \hat{p}_{\delta}(\y)   \zeta_{\delta}(\x,\y)  \frac{1}{\hat{\omega}_{\delta}(\y)} d \tau_\y =r_{\delta}(\x) \qquad \forall \ \x \in \M.
\end{equation}
Taking the derivative of \eqref{higher4} yields
\begin{equation} \label{higher6}
\begin{split}
& \frac{1}{\delta^2} \nabla^{\x}_{\M}  \int_\M p_{\delta}(\x) R_{\delta}(\x,\y)  d \mu_\y  
- \frac{1}{\delta^2} \nabla^{\x}_{\M}  \int_\M p_{\delta}(\y) R_{\delta}(\x,\y)  d \mu_\y  -
 \delta C_R \int_{\partial \M} \Delta_{ \partial \M} \hat{p}_{\delta}(\y)  ( \nabla^{\x}_{\M}   \zeta_{\delta}(\x,\y))  \frac{1}{\hat{\omega}_{\delta}(\y)} d \tau_\y  = \nabla_{\M} r_{\delta}(\x).
 \end{split}
\end{equation}
Recall the definition of $\tilde{p}_{\delta}$ and $\omega_1$ in \eqref{tildep}, we rewrite the first two terms of \eqref{higher6} and applying integration by parts to the boundary term to reformulate \eqref{higher6} as
\begin{equation} \label{higher7}
\begin{split}
\!\!\!\!\! \nabla_{\M}  (p_{\delta}(\x) \omega_1(\x))
-  \nabla_{\M}  (\tilde{p}_{\delta}(\x)  \omega_1(\x))+  \delta^3 C_R \int_{\partial \M} \nabla_{ \partial \M} \hat{p}_{\delta}(\y) \cdot \nabla_{ \partial \M}^{\y} \big( ( \nabla^{\x}_{\M}   \zeta_{\delta}(\x,\y))  \frac{1}{\hat{\omega}_{\delta}(\y)} \big) d \tau_\y 
 = \delta^2 \nabla_{\M} r_{\delta}(\x).
 \end{split}
\end{equation}
Let us separate the term $\nabla_{\M}  p_{\delta}$:
\begin{equation} \label{higher8}
\begin{split}
 \nabla_{\M}  p_{\delta}(\x) =  & \nabla_{\M}  \tilde{p}_{\delta}(\x) - \frac{\nabla_{\M} \omega_1(\x)}{\omega_1(\x)} (p_{\delta}(\x) - \tilde{p}_{\delta}(\x) )  + \delta^2 \frac{ \nabla_{\M} r_{\delta}(\x)} {\omega_1(\x)}
 + \delta^3 C_R \int_{\partial \M}    \nabla_{ \partial \M}   \hat{p}_{\delta}(\y)  \cdot \frac{1}{\omega_1(\x) } \nabla^{\y}_{\partial \M} \frac{  \nabla^{\x}_{\M}  \zeta_{\delta} (\x, \y)   }{ \hat{\omega}_{\delta}(\y) }   d \tau_\y.
 \end{split}
\end{equation}
The definition of $\omega_1$ in \eqref{tildep} implies the existence of $C_1, C_2>0$ such that $C_1  \leq \omega_1(\x) \leq C_2 $ for all $\x \in \M$, and
\begin{equation} 
\begin{split}
|\nabla_{\M} \omega_1(\x)|= & |\int_{\M} \nabla^{\x}_{\M} R_{\delta} (\x, \y) d \mu_{\y}|
\leq \frac{1}{2\delta^2} \int_{\M}   |\x-\y| \underline{R}_{\delta} (\x, \y) d \mu_{\y}
\leq C \delta^{-1}
\int_{\M}  \underline{R}_{\delta} (\x, \y)  d \mu_{\y} \leq C \delta^{-1},
\end{split}
\end{equation}
 Using the above bound of $\omega_1$, we derive from \eqref{higher8} that
\begin{equation}  \label{higher10}
\begin{split}
\!\!\!\!\!\!\!\!\!\!  | \nabla_{\M}  p_{\delta}(\x) | \leq C ( | \nabla_{\M}  \tilde{p}_{\delta}(\x)| +\frac{1}{\delta} |p_{\delta}(\x) - \tilde{p}_{\delta}(\x)|  + \delta^2 | \nabla_{\M} r_{\delta}(\x)|
 + \delta^3  \int_{\partial \M}   \big| \nabla_{ \partial \M}   \hat{p}_{\delta}(\y) \big| \  \big| \nabla^{\y}_{\partial \M} \frac{  \nabla^{\x}_{\M}  \zeta_{\delta} (\x, \y)   }{ \hat{\omega}_{\delta}(\y) } \big|  d \tau_\y). 
 \end{split}
\end{equation}
Our goal in this lemma is to establish an $L^2$ bound for $\nabla_{\M}  p_{\delta}$. The main difficulty here is the control of the last term of \eqref{higher10}. In fact, 
\begin{equation}  \label{higher10b}
\begin{split}
& \big| \nabla^{\y}_{\partial \M} \frac{  \nabla^{\x}_{\M}  \zeta_{\delta} (\x, \y)    }{ \hat{\omega}_{\delta}(\y) } \big| 
=2\delta^2 \Big| \nabla^{\y}_{\partial \M} \frac{  \nabla^{\x}_{\M} ( \n(\y) \cdot \nabla^{\x}_{\M} \overset{=}{R}_{\delta} (\x, \y)  )   }{  \int_{\partial \M}  \n(\x) \cdot \n(\y) \overset{=}{R}_{\delta}(\x,\y) d\mu_{\x} } \Big|
\leq C\delta \big| \nabla^{\y}_{\partial \M}  \nabla^{\x}_{\M}    \nabla^{\x}_{\M}  (\n(\y)  \overset{=}{R}_{\delta} (\x, \y) ) \big|
\\
& + C\delta^2 \big| \nabla^{\y}_{\partial \M} \int_{\partial \M}  \n(\x) \cdot \n(\y) \overset{=}{R}_{\delta}(\x,\y) d\mu_{\x}   \big| \ \big| \nabla^{\x}_{\M} \nabla^{\x}_{\M} (\n(\y)  \overset{=}{R}_{\delta} (\x, \y)  ) \big| 
\leq \frac{C}{\delta^2}  ( \overset{=}{R}_{\delta} (\x, \y) +\bar{R}_{\delta}(\x,\y) + {R}_{\delta} (\x, \y)+ \underline{R}_{\delta} (\x, \y)).
 \end{split}
\end{equation}
By the same argument as the control on the second term of \eqref{higher15}, we can derive
\begin{equation}  \label{higher11}
\begin{split}
& \left \| \delta^3  \int_{\partial \M}   \big| \nabla_{ \partial \M}   \hat{p}_{\delta}(\y) \big| \  \big| \nabla^{\y}_{\partial \M} \frac{  \nabla^{\x}_{\M}  \zeta_{\delta} (\x, \y)   }{ \hat{\omega}_{\delta}(\y) } \big|  d \tau_\y \right \|^2_{L^2_{\x}(\M)}
\\ \leq &
C\delta^2 \left \|   \int_{\partial \M}   \big| \nabla_{ \partial \M}   \hat{p}_{\delta}(\y) \big|  ( \overset{=}{R}_{\delta} (\x, \y) +\bar{R}_{\delta}(\x,\y) + {R}_{\delta} (\x, \y)+ \underline{R}_{\delta} (\x, \y)) d \tau_\y  \right \|^2_{L^2_{\x}(\M)}
\leq C \delta  \left \| \nabla_{ \partial \M}   \hat{p}_{\delta} \right \|^2_{L^2(\partial \M)},
 \end{split}
\end{equation}
which is controlled above by $C B_{\delta}[p_{\delta}, p_{\delta}]$ according to \eqref{bilinear5}. We then return to \eqref{higher10}. Applying \eqref{bilinear6}, \eqref{bilinear9} and \eqref{higher11} to the first, second and fourth term of \eqref{higher10} yields
\begin{equation}  \label{higher12}
\begin{split}
\left \| \nabla_{\M}  p_{\delta} \right \|^2_{L^2(\M)} \leq C ( B_{\delta}[p_{\delta}, p_{\delta}]+ \delta^4 \left \| \nabla_{\M}  r_{\delta} \right \|^2_{L^2(\M)}).
 \end{split}
\end{equation}
Furthermore, incorporating the $L^2$-bound on $p_\delta$ from \eqref{lax3}, we obtain the final estimate:
\begin{equation} \label{higher14}
\left \|   p_{\delta}  \right \|^2_{H^1(\M)} \leq C_1 ( B_{\delta}[p_{\delta}, p_{\delta}] +  \delta^4   \left \| \nabla_{\M} r_{\delta}  \right \|_{L^2(\M)} )
\end{equation}
where $C_1 > 0$ is a constant depending only on the geometry of $\M$.

Returning to the assumption of Lemma \ref{higher1}, let us seek an upper bound for the bilinear form $B_{\delta}[p_{\delta}, p_{\delta}]$. The equation \eqref{higher2} and \eqref{higher3} gives
\begin{equation}  \label{higher13b}
B_{\delta}[p_{\delta}, p_{\delta}]=(r_{\delta}, p_\delta) \leq C \left \| r_{1\delta} \right \|_{H^m(\M)}  (\left \| p_{\delta} \right \|_{H^1(\M)} + \left \| \bar{p}_{\delta} \right \|_{H^1(\M)}  + \left \| \overset{=}{p}_{\delta} \right \|_{H^1(\M)} ),
\end{equation}
where $ \bar{p}_{\delta} $ and $\overset{=}{p}_{\delta}$ are defined in \eqref{omega2} and \eqref{omega3}, respectively, with $q_{\delta}$ replaced by $p_{\delta}$. 

We now estimate $\left \| \bar{p}_{\delta} \right \|_{H^1(\M)}$. Substituting $R_{\delta}$ into $\bar{R}_{\delta}$ in \eqref{bilinear6} and \eqref{bilinear9} respectively, we obtain
\begin{equation}  \label{higher17}
\left \| \nabla \bar{p}_{\delta} \right \|_{L^2(\M)}^2 +  \left \| \bar{p}_{\delta}-p_{\delta} \right \|^2_{L^2(\M)} \leq \frac{1}{2\delta^2} \int_\M \int_\M (p_{\delta}(\x) -p_{\delta}(\y) )^2  \bar{R}_{\delta}(\x,\y)  d \mu_\y d \mu_\x.
\end{equation}
Moreover, analogous to \eqref{bilinear9}, [\cite{Yjcms1} Eq. (5.15) (5.21)] demonstrate that:
\begin{equation}  \label{higher18}
 \int_\M \int_\M (p_{\delta}(\x) -p_{\delta}(\y) )^2  \bar{R}_{\delta}(\x,\y)  d \mu_\y d \mu_\x  \leq C \int_\M \int_\M (p_{\delta}(\x) -p_{\delta}(\y) )^2  {R}_{\delta}(\x,\y)  d \mu_\y d \mu_\x .
\end{equation}
Combining the bounds from \eqref{higher17} \eqref{higher18} and \eqref{lax3}, we obtain
\begin{equation} \label{higher141}
\left \| \bar{p}_{\delta} \right \|^2_{H^1(\M)} \leq C_2 \frac{1}{\delta^2} \int_{\M}  \int_{\M} (p_{\delta} (\y)-p_{\delta}(\x))^2 R_{\delta}(\x,\y)  d \mu_\y d{\mu}_{\x} \leq C_2 B_{\delta}[p_{\delta}, p_{\delta}]
\end{equation}
for some constant $C_2>0$. Following the same analysis applied to $\bar{p}_{\delta}$, we derive
\begin{equation} \label{higher142}
\left \| \overset{=}{p}_{\delta} \right \|^2_{H^1(\M)} \leq C_3 B_{\delta}[p_{\delta}, p_{\delta}],
\end{equation}
where $C_3>0$ is another uniform constant. Incorporating these estimates for $\bar{p}_{\delta}, \overset{=}{p}_{\delta}$,and the $p_{\delta}$-bound from \eqref{higher14} into \eqref{higher13b}, we conclude
\begin{equation}  \label{higher13c}
\begin{split}
& \!\!\!\!\!\!\!\!\!\! B_{\delta}[p_{\delta}, p_{\delta}] \leq C \left \| r_{1\delta} \right \|_{H^k(\M)}  (\left \| p_{\delta} \right \|_{H^1(\M)} + \left \| \bar{p}_{\delta} \right \|_{H^1(\M)}  + \left \| \overset{=}{p}_{\delta} \right \|_{H^1(\M)} ) \\
\leq & (C^2C_1  \left \| r_{1\delta} \right \|^2_{H^k(\M)} + \frac{1}{4C_1} \left \| p_{\delta} \right \|^2_{H^1(\M)})  
+   (C^2C_2  \left \| r_{1\delta} \right \|^2_{H^k(\M)} + \frac{1}{4C_2} \left \| \bar{p}_{\delta} \right \|^2_{H^1(\M)})  
+  (C^2C_3  \left \| r_{1\delta} \right \|^2_{H^k(\M)} + \frac{1}{4C_3} \left \| \overset{=}{p}_{\delta} \right \|^2_{H^1(\M)}) \\
\leq & C^2( C_1 +C_2+C_3) \left \| r_{1\delta} \right \|^2_{H^k(\M)} + \frac{3}{4} B_{\delta}[p_{\delta}, p_{\delta}] +\frac{1}{4} \delta^4   \left \| \nabla_{\M} r_{\delta}  \right \|^2_{L^2(\M)},
\end{split}
\end{equation}
Here, $C_1, C_2, C_3$ denote the constants appearing in \eqref{higher14}, \eqref{higher141}, and \eqref{higher142}, respectively.
By subtracting $\frac{3}{4}B_{\delta}[p_{\delta}, p_{\delta}]$ from both sides of \eqref{higher13c}, we derive:
\begin{equation}  \label{higher13}
B_{\delta}[p_{\delta}, p_{\delta}] \leq C (\left \| r_{1\delta} \right \|^2_{H^k(\M)} + \delta^4   \left \| \nabla_{\M} r_{\delta}  \right \|^2_{L^2(\M)}).
\end{equation}
Applying the bound from \eqref{higher14} to \eqref{higher13} yields:
\begin{equation}
\left \|    p_{\delta}  \right \|_{H^1(\M)}  \leq C  ( \left \| r_{1\delta} \right \|_{H^k(\M)}+ \delta^2   \left \| \nabla_{\M} r_{\delta}  \right \|_{L^2(\M)}).
\end{equation}
This completes the proof of lemma.
\end{proof}

\section{Proof of Equality \eqref{mm10}} \label{appen0}
The purpose of this section is to present
$\Delta_{\M} u(\y)= \Delta_{\partial \M} u(\y)+  \nabla_{\M}^2 u(\y) (\n,\n)$ for each $\y \in \partial \M$.

 For each fixed $ \y_0 \in \partial\mathcal{M}$, there exists a neighborhood $ U \subset \mathcal{M} $ of $ \y_0 $ such that $ U$ admits a boundary normal coordinate system
$$
\y = \Phi(\theta_1, \dots, \theta_{m-1}, t),
$$
where:
\begin{enumerate}
\item
$ \Phi(\theta_1, \dots, \theta_{m-1}, 0) $ parametrizes $ \partial \M \cap U$;
\item $t$ is the signed geodesic distance from $ \y$ to $ \partial \M$, measured along the unique minimizing geodesic in $\M$ that intersects $\partial \M$ orthogonally;
\item The curves $t \mapsto \Phi(\bm{\theta}, t)$ are unit-speed geodesics normal to $\partial \M$, where $\bm{\theta}=(\theta_1, \dots, \theta_{m-1})$;
\item In these coordinates, the metric takes the form  
$$
g = \begin{pmatrix}
g_{\partial}(\bm{\theta},t) & 0 \\
0 & 1
\end{pmatrix}
$$
where $g_{\partial}(\bm{\theta}, t)$ is the induced metric on the level set $\{t = \text{constant}\}$.
\end{enumerate}
Recall the definition of $div_{\M}$ for any vector field $F: \M \to \mathcal{T}^{\y}_{\M}$ at the point $\y=\Phi(\bm{\theta}, t)=\Phi(\theta_1, ..., \theta_{m-1}, \theta_m)$:
\begin{equation} \label{divF}
div_{\M}(F)=\frac{1}{\sqrt{det \ g}} \sum \limits_{i,j=1}^m \frac{\partial}{\partial \theta_i} (\sqrt{det \ g} \ g^{ij} F^k(\y) \frac{\partial \Phi^k}{\partial \theta_j}),
\end{equation}
where $g^{ij}= \{ g^{-1} \}_{ij}, \ i,j=1,2,...m.$
Using the boundary normal coordinates, we calculate the following term at $\y_0=\Phi(\bm{\theta},0)$:
\begin{equation} \label{A.2}
\begin{split}
\Delta_{\M} u= & \text{div}_{\M} (\nabla_{\M} u) = \frac{1}{\sqrt{\det g}} \left[ \sum_{i,j=1}^{m-1}  \frac{\partial}{\partial \theta_i}(  g^{ij} \sqrt{\det g} ( \nabla_{\M} u \cdot \frac{\partial \Phi}{\partial \theta_j} )) + \frac{\partial}{\partial t}( g^{tt} \sqrt{\det g} ( \nabla_{\M} u \cdot \frac{\partial \Phi}{\partial t}  )) \right] \\
= & \frac{1}{\sqrt{\det g_{\partial }}}   \sum_{i,j=1}^{m-1}  \frac{\partial}{\partial \theta_i}(  g_{\partial}^{ij} \sqrt{\det g_{\partial} } ( \nabla_{ \M} u \cdot \frac{\partial \Phi}{\partial \theta_j} )) + \frac{1}{\sqrt{\det g_{\partial }}} \frac{\partial}{\partial t}( \sqrt{\det g_{\partial}} ( \nabla_{\M} u \cdot \frac{\partial \Phi}{\partial t}  )) \\
 = & \text{div}_{\partial \M} (\nabla_{\M} u |_{\partial \M} ) + (\frac{1}{\sqrt{\det g_{\partial}}}  \frac{\partial }{\partial t} \sqrt{\det g_{\partial}}  ) ( \nabla_{\M} u \cdot \frac{\partial \Phi}{\partial t}  )
 +   \frac{\partial }{\partial t} ( \nabla_{\M} u \cdot \frac{\partial \Phi}{\partial t}  ) ,
\end{split}
\end{equation}
where we have used the fact that $g^{tt} \equiv 1$ and $g^{it} \equiv 0, i=1,2,...,m-1$ in the second and third equal sign. 
When under the condition $\frac{\partial u}{\partial \n}=0$, we have $\nabla_{\M} u |_{\partial \M}=\nabla_{\partial \M} u$. Besides, the second term of \eqref{A.2} vanishes due to $\frac{\partial \Phi}{\partial t} =\n$. We then continue writing \eqref{A.2} as
\begin{equation}
\Delta_{\M} u= \text{div}_{\partial \M} (\nabla_{\partial \M} u)+ \frac{\partial }{\partial t} ( \nabla_{\M} u \cdot \frac{\partial \Phi}{\partial t}  ) = \Delta_{\partial \M} u+  \nabla_{\M}^2 u (\n,\n);
\end{equation}
which is exactly \eqref{mm10}.

\section{Lemmas for estimation of \eqref{e090}} \label{appenb}
This section introduces some lemmas that helps to control the truncation error $\Delta_{\partial \M}(u-\hat{u})$.

\begin{lmm}  \label{lmm2}
Let $\overset{=}{R}_{\delta}$ denote the kernel function defined in \eqref{kernel}. Then there exists a constant $C > 0$ such that for all $\y \in \partial \M$:
\begin{equation}
 |  \Delta^{\y}_{\partial \M}  \int_{ \M}  \overset{=}{R}_{\delta} (\x, \y)  d \tau_{\x} | + |  \nabla^{\y}_{\partial \M}   \int_{ \M} \overset{=}{R}_{\delta} (\x, \y)  d \tau_{\x} |  \leq C \delta.
\end{equation}
\end{lmm}

\begin{lmm}[Interior Function Estimate]  \label{lmm5}
 For any function $p \in H^2(\M)$, define $$p_1(\y)=\int_{\M} p(\x) \overset{=}{R}_{\delta} (\x, \y)   d \mu_\x, \qquad \forall \  \y \in \partial \M.$$ Then the following estimate holds 
\begin{equation} 
\left \| \nabla_{\partial \M} p_1 \right \|_{L^2(\partial \M)} + \left \| \Delta_{\partial \M} p_1 \right \|_{L^2(\partial \M)} \leq C \delta^{-\frac{1}{2}} \left \|  p \right \|_{H^2(\M)} .
\end{equation}
\end{lmm}

\begin{lmm}  \label{lmm3}
Let $\overset{=}{R}_{\delta}$ denote the kernel function defined in \eqref{kernel}. Then, 
\begin{enumerate}
\item
there exists a constant $C > 0$ such that:
\begin{equation}
| 2\delta \int_{\partial \M}  \overset{=}{R}_{\delta} (\x, \y)  d \tau_{\x} - C_R | \leq C\delta^2, \qquad \forall \ \y \in \partial \M,
\end{equation}
where $C_R$ is defined in \eqref{defome1};
\item
there exists a constant $C > 0$ such that:
\begin{equation}
 | \Delta^{\y}_{\partial \M}   \int_{\partial \M}  \overset{=}{R}_{\delta} (\x, \y)  d \tau_{\x} | + | \nabla^{\y}_{\partial \M}  \int_{\partial \M}   \overset{=}{R}_{\delta} (\x, \y)  d \tau_{\x} |  \leq C \delta, \qquad \forall \ \y \in \partial \M.
\end{equation}
\end{enumerate}
\end{lmm}

\begin{lmm}[Boundary Function Estimate]  \label{lmm4}
 For any function $p \in H^2(\partial \M)$, define $$p_2(\y)=\int_{\partial \M} p(\x)  \overset{=}{R}_{\delta} (\x, \y)   d \tau_\x, \qquad \forall \  \y \in \partial \M.$$ Then the following estimate holds 
\begin{equation} 
\left \| \nabla_{\partial \M} p_2 \right \|_{L^2(\partial \M)} + \left \| \Delta_{\partial \M} p_2 \right \|_{L^2(\partial \M)} \leq C \delta^{-\frac{3}{2}} \left \|  p \right \|_{H^2(\partial \M)} .
\end{equation}
\end{lmm}

Lemma \ref{lmm2}$-$\ref{lmm4} give a complete control on the interior and boundary weighted average functions with respect to the kernel $\overset{=}{R}_{\delta}$ up to second order derivative.

\begin{proof}[Proof of Lemma \ref{lmm2}]

Fix $\y \in \partial \M$ and use the boundary normal coordinate system $\Phi: \mathbb{R}^m \to U \subset \M$ from Appendix \ref{appen0}. We assume, without loss of generality, that the induced matrix $g$ at $\y$ is the identity. Choose $\delta$ small enough that $B_{\y}(2\delta) \cap \M \subset U$. Write $\bm{\theta}'=(\bm{\theta},0)=\Phi^{-1}(\y)$. For $\x \in U$, set $\vv'=\frac{1}{2\delta} (\Phi^{-1}(\x)- \bm{\theta}') $ so that $\x=\Phi(\bm{\theta}'+2\delta \vv')$ and $\y=\Phi(\bm{\theta}')=\Phi(\bm{\theta},0)$. Then
\begin{equation} \label{B.1}
\begin{split}
F_{\delta}(\y)= \int_{\M}  \overset{=}{R}_{\delta} (\x, \y)  d \mu_{\x} =   \int_{\M}  \overset{=}{R} (\frac{|\x-\y|^2}{4\delta^2}) d \mu_{\x} 
 = \frac{(2\delta)^m}{(4\pi \delta^2)^{\frac{m}{2}} }\int_{\mathbb{R}^m \cap \{t \geq 0\} }  \overset{=}{R}  (\frac{ | \Phi (\bm{\theta}' +2\delta \vv') - \Phi(\bm{\theta}' )|^2 } {4\delta^2} ) \sqrt{\mbox{det} \ g(\bm{\theta}'+2\delta \vv')} d  \vv'.
\end{split}
\end{equation}
Because $\Phi$ is smooth, Taylor expansion gives
\begin{equation} \label{B.2}
\begin{split}
\Phi (\bm{\theta}' +2\delta \vv') - \Phi(\bm{\theta}' )=2\delta \ \nabla \Phi(\bm{\theta'}) ( \vv')  +2\delta^2 \ \nabla^2 \Phi(\bm{\theta'})( \vv', \vv')
+ \mathcal{O}(\delta^3).
\end{split}
\end{equation}
Since $\overset{=}{R} \in C^4$ as we assumed, we obtain
\begin{equation} \label{B.3}
\begin{split}
 & \overset{=}{R}  (\frac{ | \Phi (\bm{\theta}' +2\delta \vv') - \Phi(\bm{\theta}' )|^2 } {4\delta^2} ) = 
  \overset{=}{R} \bigg( \big|    \nabla \Phi(\bm{\theta}') (\vv')  + \delta  \nabla^2 \Phi(\bm{\theta}')( \vv', \vv') \big|^2 + \mathcal{O}(\delta^2) \bigg) \\ 
& =  \overset{=}{R} \bigg(   |\vv'|^2 + 2  \delta \ \nabla \Phi(\bm{\theta}') (\vv')  \cdot \nabla^2 \Phi(\bm{\theta}')( \vv', \vv') \bigg) + \mathcal{O}(\delta^2)  
 =  \overset{=}{R} \big(  |\vv'|^2 \big)- 2 \delta \ \bar{R} \big(   |\vv'|^2 \big) \ \nabla \Phi(\bm{\theta}') (\vv')  \cdot \nabla^2 \Phi(\bm{\theta}')( \vv', \vv') + \mathcal{O}(\delta^2) ,
\end{split}
\end{equation}
where we used $|\nabla \Phi(\bm{\theta}') (\vv') | =| \vv'|$ (metric is identity at $\y$). 

For the volume element,
\begin{equation} \label{B.4}
\begin{split}
\sqrt{\mbox{det} \ g(\bm{\theta}'+2\delta \vv')} =\sqrt{\mbox{det} \ g(\bm{\theta}')} +2\delta \vv' \cdot  \frac {(\nabla  \mbox{det} \ g) (\bm{\theta}') }{ 2 \sqrt{\mbox{det} \ g(\bm{\theta}')}   }  + \mathcal{O}(\delta^2)=1+\delta \vv' \cdot (\nabla \mbox{det} \ g) (\bm{\theta}')+ \mathcal{O}(\delta^2).
\end{split}
\end{equation}
Substituting the above two estimates into \eqref{B.1} yields
\begin{equation}
\begin{split}
F_{\delta}(\y)= \int_{\M}  \overset{=}{R}_{\delta} (\x, \y)  d \mu_{\x} 
 =C_0+ \delta C_1(\y) + \mathcal{O}(\delta^2),
\end{split}
\end{equation}
with
\begin{equation}
C_0=\pi^{-\frac{m}{2}} \int_{\mathbb{R}^m \cap \{t \geq 0\} }  \overset{=}{R} \big(  |\vv'|^2 \big) d \vv',
\end{equation}
and
\begin{equation} \label{B.6}
C_1(\y) =\pi^{-\frac{m}{2}} \int_{\mathbb{R}^m \cap \{t \geq 0\} }    
 \ \overset{=}{R} \big(   |\vv'|^2 \big) \ \vv' \cdot \nabla (\mbox{det} \ g) (\bm{\theta}')
- 2 \bar{R} \big(   |\vv'|^2 \big) \ \nabla \Phi(\bm{\theta}') (\vv')  \cdot \nabla^2 \Phi(\bm{\theta}')( \vv', \vv')
 \   d \vv'.
\end{equation}
Here $\bm{\theta}'=(\bm{\theta},0)=\Phi^{-1}(\y)$. Because $\mbox{det} \ g$ and $\Phi$ are smooth, $C_1(\y)$ is independent of $\delta$ and smooth on $\y \in \partial \M$. In fact, if $\bar{R} \in {C}^k[0,1]$, we have the expansion
\begin{equation} \label{B.5}
F_{\delta}(\y)=C_0+\sum_{q=1}^k \delta^q C_q(\y)+ \mathcal{O}(\delta^{k+1}),
\end{equation}
where each $C_q: \partial \M \to \mathbb{R}, q=1,2,...,k$ is smooth and  independent of $\delta$. Our assumption on $\bar{R}$ gives $k \geq 3$. 

Recall the surface gradient $\nabla_{\partial \M}$ and Laplace-Beltrami operator $\Delta_{\partial \M}$ for any function $\bar{C}: \partial \M \to \mathbb{R}$:
\begin{equation} \label{nablaF}
\begin{split}
& \nabla_{\partial \M}  \bar{C}(\y) = \sum \limits_{i,j=1}^{m-1} g_{\partial}^{ij}(\bm{\theta}, 0) \frac{\partial \Phi(\bm{\theta},0) }{\partial \theta_i}  \frac{\partial \bar{C}(\Phi(\bm{\theta}, 0))}{\partial \theta_j}; \\
& \Delta_{\partial \M} \bar{C} (\y) = \frac{1}{\sqrt{\mbox{det} \ g_{\partial} (\bm{\theta},0)  }} \sum \limits_{l,s=1}^{m-1} \frac{\partial}{\partial \theta_l} (\sqrt{\mbox{det} \ g_{\partial} (\bm{\theta},0)  } g_{\partial}^{ls} (\bm{\theta},0)  \sum \limits_{i,j=1}^{m-1} g_{\partial}^{ij}(\bm{\theta}, 0) \frac{\partial \Phi(\bm{\theta},0) }{\partial \theta_i}  \frac{\partial \bar{C}(\Phi(\bm{\theta}, 0))}{\partial \theta_j} \cdot \frac{\partial \Phi(\bm{\theta},0) }{\partial \theta_s}),
\end{split}
\end{equation}
where $\y=\Phi(\bm{\theta,0})$. This derivatives act only on the $\y$-variable(or $\bm{\theta}$) and do not involve $\delta$.
Applying \eqref{B.5} with $k = 1$ gives
\begin{equation} \label{B.15}
\begin{split}
\nabla_{\partial \M} F_{\delta}(\y)= \delta \nabla_{\partial \M} C_1(\y) +\mathcal{O}(\delta^2), \qquad
\Delta_{\partial \M} F_{\delta}(\y)= \delta \Delta_{\partial \M} C_1(\y) +\mathcal{O}(\delta^2).
\end{split}
\end{equation}
\eqref{B.6} indicates that $C_1(\y)=C_1(\Phi(\bm{\theta},0))$ is smooth on $\bm{\theta}$ on the compact manifold $\partial \M$. According to \eqref{nablaF}, its derivatives on $\y$ are also uniformly bounded.
 Lemma \ref{lmm2} thereby hold.

\end{proof}

\begin{proof}[Proof of Lemma \ref{lmm5}]
Continue with the proof of Lemma \ref{lmm2}, for each $\y \in \partial \M$, we follow the same reasoning as \eqref{B.1} to compute
\begin{equation} \label{B.7}
\begin{split}
p_1(\y)=p_1(\Phi(\bm{\theta}'))=
\frac{(2\delta)^m}{(4\pi \delta^2)^{\frac{m}{2}} }\int_{\mathbb{R}^m \cap \{t \geq 0\} }  p(\Phi (\bm{\theta}' +2\delta \vv')) \overset{=}{R}  (\frac{ | \Phi (\bm{\theta}' +2\delta \vv') - \Phi(\bm{\theta}' )|^2 } {4\delta^2} ) \sqrt{\mbox{det} \ g(\bm{\theta}'+2\delta \vv')} d  \vv'.
\end{split}
\end{equation}
We rewrite it as
\begin{equation}
\begin{split}
p_1(\y)=\pi^{-\frac{m}{2}} \int_{\mathbb{R}^m \cap \{t \geq 0\} } \beta(\bm{\theta}', \vv', \delta) d \vv',
\end{split}
\end{equation}
where $\beta$ is the integrand of \eqref{B.7}. According to \eqref{nablaF}, the function $\nabla_{\partial \M} p'_1(\y)$ and $\Delta_{\partial \M} p'_1(\y)$ can be formulated as the derivatives of $p'_1(\Phi(\bm{\theta},0))$ on $\bm{\theta}$ multiplied by a $C^{\infty}$ tensor. 
Hence our task now is to bound the $\bm{\theta}$ derivative of each term that composes $\beta$. 
Recall the decomposition \eqref{B.3}, we are able to deduce
\begin{equation} \label{B.8}
\begin{split}
\overset{=}{R}  (\frac{ | \Phi (\bm{\theta}' +2\delta \vv') - \Phi(\bm{\theta}' )|^2 } {4\delta^2} ) + \Big| \frac{\partial}{\partial \theta_i} \overset{=}{R}  (\frac{ | \Phi (\bm{\theta}' +2\delta \vv') - \Phi(\bm{\theta}' )|^2 } {4\delta^2} ) \Big| +  \Big| \frac{\partial^2}{\partial \theta_i \partial \theta_j} \overset{=}{R}  (\frac{ | \Phi (\bm{\theta}' +2\delta \vv') - \Phi(\bm{\theta}' )|^2 } {4\delta^2} ) \Big| \leq C
\end{split}
\end{equation}
for any $1 \leq i,j \leq m-1$.
By a similar argument, the Taylor expansion \eqref{B.4} implies
\begin{equation} \label{B.9}
\begin{split}
\sqrt{\mbox{det} \ g(\bm{\theta}'+2\delta \vv')} + \Big| \frac{\partial}{\partial \theta_i} \sqrt{\mbox{det} \ g(\bm{\theta}'+2\delta \vv')}   \Big|  +  \Big| \frac{\partial^2}{\partial \theta_i \partial \theta_j } \sqrt{\mbox{det} \ g(\bm{\theta}'+2\delta \vv')}   \Big| \leq C, \quad \forall \ 1 \leq i,j \leq m-1.
\end{split}
\end{equation}
For the last term that composes $\beta$, we apply the formula on derivatives of composite functions to obtain
\begin{equation} \label{B.12}
\begin{split}
\big|  \frac{\partial}{\partial \theta_i}  p(\Phi (\bm{\theta}' +2\delta \vv')) \big| + \big| \frac{\partial^2}{\partial \theta_i \partial \theta_j } \big( p(\Phi (\bm{\theta}' +2\delta \vv')) \big| \leq C \sum_{q=1}^2 \big| ( \nabla^q_{\M} p) (\Phi (\bm{\theta}' +2\delta \vv')) \big| , \qquad \forall \ 1 \leq i,j \leq m-1.
\end{split}
\end{equation}
Finally, we combine all the above estimate \eqref{B.8}-\eqref{B.12} and return to \eqref{B.7} to obtain
\begin{equation} \label{B.13}
\begin{split}
& |\nabla_{\partial \M} p_1(\y)|+  |\Delta_{\partial \M} p_1(\y)|
 \leq C \int_{\mathbb{R}^m \cap \{ t \geq 0\} }  \sum_{i=1}^{m-1} \big| \frac{\partial}{\partial \theta_i}  \beta(\bm{\theta}', \vv', \delta) \big| +\sum_{i,j=1}^{m-1}  \big|\frac{\partial^2}{\partial \theta_i \partial \theta_j}  \beta(\bm{\theta}', \vv', \delta) \big|  d \vv' \\
 & =  C \int_{ \substack{ \mathbb{R}^m \cap \{ t \geq 0\} \\ | \Phi (\bm{\theta}' +2\delta \vv') - \Phi(\bm{\theta}' )|  \leq 2\delta }  }  \sum_{i=1}^{m-1} \big| \frac{\partial}{\partial \theta_i}  \beta(\bm{\theta}', \vv', \delta) \big| +\sum_{i,j=1}^{m-1}  \big|\frac{\partial^2}{\partial \theta_i \partial \theta_j}  \beta(\bm{\theta}', \vv', \delta) \big|  d \vv'  \\ 
& \leq C \int_{ \substack{ \mathbb{R}^m \cap \{ t \geq 0\} \\ | \Phi (\bm{\theta}' +2\delta \vv') - \Phi(\bm{\theta}' )|  \leq 2\delta }  } \sum_{q=0}^2   \big| ( \nabla^q_{\M} p) (\Phi (\bm{\theta}' +2\delta  \vv')) \big|  d \vv' \leq C (2\delta)^{-m} \int_{\M \cap  B_{\y}(2\delta) } (|p(\x)| + |\nabla_{\M} p(\x)| +|\nabla^2_{\M} p(\x)| ) d \mu_{\x},
\end{split}
\end{equation}
where in the equality sign of \eqref{B.13} we have used the fact that $\beta(\bm{\theta}', \vv', \delta) \equiv 0$ when $| \Phi (\bm{\theta}' +2\delta \vv') - \Phi(\bm{\theta}' )| > 2\delta.$ By the randomness of $\y$, \eqref{B.13} hold for all $\y \in \partial \M$. Using Cauchy-Schwarz inequality, this yields:
\begin{equation} 
\begin{split}
& \left \| \nabla_{\partial \M} p_1 \right \|^2_{L^2(\partial \M)} + \left \| \Delta_{\partial \M} p_1 \right \|^2_{L^2(\partial \M)} \\
 \leq & C \delta^{-2m} \int_{\partial \M} (\int_{\M \cap  B_{\y}(2\delta) }  1 d \mu_{\x} \int_{\M \cap  B_{\y}(2\delta) } (|p(\x)|^2 + |\nabla_{\M} p(\x)|^2 +|\nabla^2_{\M} p(\x)|^2 ) d \mu_{\x}) d \tau_{\y} \\
 \leq & C \delta^{-m} \int_{\M} \int_{\partial \M \cap B_{\x}(2\delta)} d \tau_\y  (|p(\x)|^2 + |\nabla_{\M} p(\x)|^2 +|\nabla^2_{\M} p(\x)|^2 ) d \mu_{\x}
 \leq C \delta^{-1} \left \|  p \right \|^2_{H^2(\M)}.
\end{split}
\end{equation}
We then complete the proof by taking the square root of both sides.
\end{proof}

\begin{proof}[Proof of Lemma \ref{lmm3}]
For brevity, we use the same boundary normal coordinate $\Phi$ as in Lemma \ref{lmm2}. This time, both $\bm{\theta}'$ and $\vv'$ have no normal components. As a consequence,
\begin{equation} 
\begin{split}
G_{\delta}(\y)= \int_{\partial \M}  \overset{=}{R}_{\delta} (\x, \y)  d \tau_{\x} =   \int_{\partial \M}  \overset{=}{R} (\frac{|\x-\y|^2}{4\delta^2}) d \tau_{\x} 
 = \frac{(2\delta)^{m-1}}{(4\pi \delta^2)^{\frac{m}{2}} }\int_{\mathbb{R}^{m-1}}   \overset{=}{R}  (\frac{ | \phi (\bm{\theta} +2\delta \vv) - \phi(\bm{\theta} )|^2 } {4\delta^2} ) \sqrt{\mbox{det} \ g_{\partial} (\bm{\theta}+2\delta \vv)} d  \vv.
\end{split}
\end{equation}
where $\phi: \mathbb{R}^{m-1} \to U \cap \partial \M$ is the boundary normal coordinate system $\Phi$ restricted on $t=0$, with $\x= \phi(\bm{\theta} +2\delta \vv), \y=\phi(\bm{\theta})$.  We then follow the same calculation as $F_{\delta}(\y)$ in Lemma \ref{lmm2} to obtain
\begin{equation}
\begin{split}
G_{\delta}(\y)= \int_{\M}  \overset{=}{R}_{\delta} (\x, \y)  d \tau_{\x} 
 =\frac{1}{2\delta} (D_0+ \delta D_1(\y) + \mathcal{O}(\delta^2)),
\end{split}
\end{equation}
with
\begin{equation}
D_0=\pi^{-\frac{m}{2}} \int_{\mathbb{R}^{m-1}} \overset{=}{R} \big(  |\vv|^2 \big) d \vv,
\end{equation}
and
\begin{equation} \label{B.14}
D_1(\y) =\pi^{-\frac{m}{2}} \int_{\mathbb{R}^{m-1}}    
 \ \overset{=}{R} \big(   |\vv|^2 \big) \ \vv \cdot \nabla (\mbox{det} \ g_{\partial} ) (\bm{\theta})
- 2 \bar{R} \big(   |\vv|^2 \big) \ \nabla \phi(\bm{\theta}) (\vv)  \cdot \nabla^2 \phi(\bm{\theta})( \vv, \vv)
 \   d \vv.
\end{equation}
It can be observed that $D_0=C_R$, and $D_1(\y) \equiv 0$ as the integrand in \eqref{B.14} is odd symmetric in ${\mathbb{R}^{m-1}}$. This indicates
\begin{equation}
2\delta G_{\delta}(\y)= D_0+\delta D_1(\y) +\mathcal{O}(\delta^2) =C_R+ \mathcal{O}(\delta^2) .
\end{equation}
Besides, utilizing the same argument as \eqref{nablaF} and \eqref{B.15} in Lemma \ref{lmm2}, we are able to conclude
\begin{equation}
\begin{split}
\nabla_{\partial \M}  G_{\delta}(\y)= \mathcal{O}(\delta), \qquad
\Delta_{\partial \M}  G_{\delta}(\y)= \mathcal{O}(\delta).
\end{split}
\end{equation}
We then complete the proof of Lemma \ref{lmm3}.

\end{proof}

We omit the proof of Lemma \ref{lmm4} as the proof is almost identical to Lemma \ref{lmm5} restricted on $\partial \M$.

In addition to Lemma \ref{lmm2}$-$\ref{lmm4}, we have the following corollary:
\begin{crllr} \label{lmm6}
Let $\n(\y)$ be the outward unit normal of $\partial \M$ at $\y$, then for $p \in H^2(\partial \M)$, we have
\begin{equation} \label{B.17}
\begin{split}
& | \nabla_{\partial \M}^{\y} \int_{\partial \M} (1-\n(\x)  \cdot \n(\y)) \overset{=}{R}_{\delta} (\x, \y)   d \tau_\x) |
+ | \Delta_{\partial \M}^{\y} \int_{\partial \M} (1-\n(\x)  \cdot \n(\y)) \overset{=}{R}_{\delta} (\x, \y)   d \tau_\x) |
 \leq C\delta,
\\
& \left \| \int_{\partial \M} p(\x) (1-\n(\x)  \cdot \n(\y)) \overset{=}{R}_{\delta} (\x, \y)   d \tau_\x \right \|_{H_{\y}^2(\partial \M)}  
\leq C\delta^{\frac{1}{2}}  \left \| p \right \|_{H^2(\partial \M)}.
 \end{split}
\end{equation}

\end{crllr}
\begin{proof}
We denote the coordinate $\phi_t(\bm{\theta})=\frac{\partial \Phi}{\partial t}(\bm{\theta,0})$ for $\bm{\theta} \in \mathbb{R}^{m-1}$.  For any $\x, \y \in \partial \M$ with $|\x-\y|<2\delta$, we have the following Taylor expansion:
\begin{equation} \label{B.16}
\begin{split}
\n(\y) \cdot ( \n(\y)-\n(\x)) =  \phi_t (\bm{\theta}) \cdot 
(    \phi_t (\bm{\theta}) -  \phi_t (\bm{\theta}+2\delta \vv )) 
=2\delta  \phi_t (\bm{\theta}) \cdot \nabla \phi_t (\bm{\theta})    (\vv) + 
2\delta^2  \phi_t (\bm{\theta}) \cdot \nabla^2 \phi_t (\bm{\theta}) ( \vv, \vv)+ \mathcal{O}(\delta^3),
\end{split}
\end{equation}
The first term of \eqref{B.16} vanishes since $ \phi_t (\bm{\theta}) \cdot  \phi_t (\bm{\theta})  \equiv 1$ for all $\bm{\theta}$. For the other terms in \eqref{B.16}, the derivatives on $\bm{\theta}$ does not affect its order of $\delta$. Following the same proof as in Lemma \ref{lmm2} and \ref{lmm5}, we are able to obtain \eqref{B.17}.

\end{proof}

\section{Control of $\tilde{f}_{\delta}$} \label{appenf}
In this section, we aim to present $\tilde{f}_{\delta}$ that defined in \eqref{preli} satisfies the estimate \eqref{vanish1}. Apart from the proofs in the previous lemmas, we apply local Euclid coordinates instead of the boundary normal coordinates.

For each boundary point $\y \in \partial \mathcal{M}$, we define a normal curve
\[
\gamma_{\y} = \{ \x \in \mathcal{M} \mid d(\x, \y) = d(\x, \partial \mathcal{M}) < 2\delta \}.
\]
In other words, $\gamma_{\y}$ consists of all points in the $2\delta$-neighborhood of the boundary for which $\y$ is the closest boundary point. Geometrically, $\gamma_{\y}$ is contained in the intersection of $\mathcal{M}$ with the affine plane orthogonal to the tangent space $\mathcal{T}_{\partial \mathcal{M}}^{\y}$, and it approximates the inward normal line $ \y - s \mathbf{n}(\y) $ for $s \in [0, 2\delta]$.

These curves are disjoint for distinct $\y \in \partial \mathcal{M}$, and their union forms the entire boundary layer:
\[
\mathcal{M}_{\delta} = \{ \x \in \mathcal{M} \mid d(\x, \partial \mathcal{M}) < 2\delta \} = \bigcup_{\y \in \partial \mathcal{M}} \gamma_{\y}.
\]

We can parameterize $\mathcal{M}_{\delta}$ using coordinates $(\y, s)$, where $\y \in \partial \mathcal{M}$ and $s$ is the arc length along $\gamma_{\y}$ from $\y$ to $\x$, denoted $\x = \gamma_{\y}(s)$ for $s \in [0, L(\y)]$. The function $L(\y)$ gives the length of $\gamma_{\y}$. In these coordinates, the volume form on $\mathcal{M}_{\delta}$ becomes
\begin{equation}
d\mu_{\x} = J(\y, s)  ds  d\tau_{\y},
\end{equation}
where $d\tau_{\y}$ is the volume form on $\partial \mathcal{M}$, and $J(\y, s)$ is the Jacobian accounting for the variation of the normal curves. Using the same argument as we calculate det $g$ in Lemma \ref{lmm2}, we find $J(\y, s) = 1 + \mathcal{O}(\delta^2)$ and $L(\y) = 2\delta (1 + \mathcal{O}(\delta^2))$.

We now return to estimating $\tilde{f}_{\delta}$ defined in \eqref{preli}. Since $ \int_{\M} f(\x) d \mu_{\x}=0$, we decompose $\tilde{f}_{\delta}$ as follows:
\begin{equation} \label{trun15}
\begin{split}
|\tilde{f}_{\delta}| = & \big| \int_{\M}  \int_{\M} f(\x) \bar{R}_{\delta}(\x,\z) d \mu_{\z} d \mu_{\x} 
-\int_{\M}  \int_{\partial \M} f(\z) (\x-\z) \cdot \n(\z) \bar{R}_{\delta}(\x,\z) d \tau_{\z} d \mu_{\x} 
- \int_{\M} f(\x) d \mu_{\x} \int_{|\vv'|<2\delta}  C_{\delta} \bar{R} ( \frac{ |\vv'|^2}{4\delta^2} )  d \vv' \big|  \\
 \leq &
 \big| \int_{\M \setminus \M_{\delta}}  f(\x) (\int_{\M} \bar{R}_{\delta}(\x,\z) d \mu_{\z} - \int_{|\vv'|<2\delta}  C_{\delta} \bar{R} ( \frac{ |\vv'|^2}{4\delta^2} )  d \vv')d \mu_{\x} \big| \\
 & +\int_{\partial \M} \int_0^{L(\y)} |f(\gamma_{\y}(s))-f(\y)| J(\y,s) ( \int_{|\vv'|<2\delta}  C_{\delta} \bar{R} ( \frac{ |\vv'|^2}{4\delta^2} )  d \vv'+ \int_{\M} \bar{R}_{\delta}(\gamma_{\y}(s),\z) d \mu_{\z} ) ds d \tau_{\y} \\
 & +\int_{\partial \M} \int_0^{L(\y)} J(\y,s) \int_{\partial \M} | (f(\z)-f(\y))(\gamma_{\y}(s)-\z) \cdot \n(\z) | \bar{R}_{\delta}(\gamma_{\y}(s),\z) d \tau_{\z} ds d \tau_{\y} \\
 & + \big| \int_{\partial \M} f(\y) \Big( \int_0^{L(\y)} J(\y,s)  \int_{\M} \bar{R}_{\delta}(\gamma_{\y}(s),\z) d \mu_{\z} ds   
  -\int_0^{L(\y)} J(\y,s) ds \int_{|\vv'|<2\delta}  C_{\delta} \bar{R} ( \frac{ |\vv'|^2}{4\delta^2} )  d \vv' \\
   & -\int_0^{L(\y)} J(\y,s) \int_{\partial \M} (\gamma_{\y}(s)-\z) \cdot \n(\z)  \bar{R}_{\delta}(\gamma_{\y}(s),\z) d \tau_{\z}ds
  \Big) d \tau_{\y} \big|.
\end{split}
\end{equation}
where $C_{\delta}=(4\pi\delta^2)^{-m/2}$ as mentioned in \eqref{base}.
Following a derivation analogous to Lemma \ref{lmm2}, we establish the estimate:
\begin{equation} 
\int_{\M} \bar{R}_{\delta}(\x,\z) d \mu_{\z} =  \int_{|\vv'|<2\delta, v_m<s}  C_{\delta} \bar{R} ( \frac{ |\vv'|^2}{4\delta^2} )  d \vv' +\mathcal{O}(\delta), \qquad  \forall  \ \s \in (0,L(\y)), \ \x= \gamma_{\y}(s).
\end{equation}
This yields
\begin{equation} \label{trun14}
\begin{split}
& \int_0^{L(\y)} J(\y,s) ds \int_{|\vv'|<2\delta}  C_{\delta} \bar{R} ( \frac{ |\vv'|^2}{4\delta^2} )  d \vv' 
- \int_0^{L(\y)} J(\y,s)  \int_{\M} \bar{R}_{\delta}(\gamma_{\y}(s),\z) d \mu_{\z} ds   \\
& = \int_0^{L(\y)}  \int_{|\vv'|<2\delta}  C_{\delta} \bar{R} ( \frac{ |\vv'|^2}{4\delta^2} )  d \vv' ds 
- \int_0^{L(\y)}  \int_{|\vv'|<2\delta, v_m<s}  C_{\delta} \bar{R} ( \frac{ |\vv'|^2}{4\delta^2} ) d \vv' ds
 +\mathcal{O}(\delta^2) \\
& =\int_0^{2\delta} \int_{|\vv'|<2\delta, v_m>s}  C_{\delta} \bar{R} ( \frac{ |\vv'|^2}{4\delta^2} ) d \vv' ds +\mathcal{O}(\delta^2)  
= \int_0^{2\delta}   s \int_{|\vv'|<2\delta, v_m=s}  C_{\delta} \bar{R} ( \frac{ |\vv|^2}{4\delta^2} )  d \vv  ds
+\mathcal{O}(\delta^2) \\
& = \int_0^{2\delta}   s \int_{\partial \M}  \bar{R}_{\delta}(\gamma_{\y}(s),\z) d \tau_{\z} ds
+\mathcal{O}(\delta^2) 
= -\int_0^{L(\y)} J(\y,s) \int_{\partial \M} (\gamma_{\y}(s)-\z) \cdot \n(\z)  \bar{R}_{\delta}(\gamma_{\y}(s),\z) d \tau_{\z} ds + \mathcal{O}(\delta^2) .
\end{split}
\end{equation}
Inequality \eqref{trun14} provides a bound for the final term in \eqref{trun15}. Proceeding to analyze the remaining terms in \eqref{trun15} and applying the method from Lemma \ref{lmm2}, we ultimately find:
\begin{equation}
\begin{split}
|\tilde{f}_{\delta}|
\leq & C ( \delta^2 \int_{\M \setminus \M_{\delta}}  |f(\x)| d \mu_{\x} 
 + \delta \int_{\partial \M} \int_0^{L(\y)} | \nabla_{\M} f(\gamma_{\y}(s)) |  ds d \tau_{\y}  \\
&  + \int_{\partial \M}    \int_{\partial \M \cap \mathcal{B}_{2\delta}(\y) } \delta^{-m+1} |\delta^2   \nabla_{\partial \M} f(\z)|   d \tau_{\z} d \tau_{\y} 
  +\delta^2  \int_{\partial \M} |f(\y)| d \tau_{\y} )  \\
  \leq & C(\delta^2 \int_{\M} |f(\x)| d \mu_{\x}+\delta ( \int_{\partial \M} L(\y) |\nabla_{\M} f(\y)| d \tau_{\y} + \delta^2 \int_{\M_{\delta}} | \nabla_{\M}^2 f(\x)| d \mu_{\x} ) \\
 & + \delta^2 \int_{\partial \M}  | \nabla_{\partial \M} f(\z)|   d \tau_{\z} + \delta^2  \int_{\partial \M} |f(\y)| d \tau_{\y} )
 \leq C\delta^2 \left \| f \right \|_{W^{2,1}(\M)} \leq  C\delta^2 \left \| f \right \|_{H^2(\M)}.
  \end{split}
\end{equation}
This is exactly the desired bound for $\tilde{f}_{\delta}$.

\end{appendix}

\bibliographystyle{abbrv}

\bibliography{reference}

\end{document}